\documentclass[10pt, a4paper]{article}

\usepackage{indentfirst}
\usepackage{amsmath, amssymb, amsfonts, amsthm, graphicx}
\usepackage{geometry}
\usepackage{setspace}

\usepackage{algorithm}        
\usepackage{algpseudocode}   
\makeatletter
\renewcommand{\ALG@beginalgorithmic}{\small} 
\makeatother

\usepackage{caption}
\usepackage{subfig}
\usepackage[section]{placeins}
\usepackage{float}
\usepackage{epstopdf}
\usepackage{cite}
\usepackage{xcolor} 

\newtheorem{theorem}{Theorem}[section]
\newtheorem{definition}{Definition}[section]

\newtheorem{lemma}[theorem]{Lemma}
\newtheorem{claim}{Claim}
\newtheorem{corollary}[theorem]{Corollary}

\newtheorem{remark}[theorem]{Remark}

\newtheorem{observation}[theorem]{Observation}

\renewcommand{\baselinestretch}{1.3}

\def\qed{\hfill \rule{4pt}{7pt}}

\RequirePackage[normalem]{ulem} 
\RequirePackage{color}\definecolor{RED}{rgb}{1,0,0}\definecolor{BLUE}{rgb}{0,0,1} 
\providecommand{\DIFaddbegin}{} 
\providecommand{\DIFaddend}{} 
\providecommand{\DIFdelbegin}{} 
\providecommand{\DIFdelend}{} 
\providecommand{\DIFaddbeginFL}{} 
\providecommand{\DIFaddendFL}{} 
\providecommand{\DIFdelbeginFL}{} 
\providecommand{\DIFdelendFL}{} 
\newcommand{\DIFscaledelfig}{0.5}
\RequirePackage{settobox} 
\RequirePackage{letltxmacro} 
\newsavebox{\DIFdelgraphicsbox} 
\newlength{\DIFdelgraphicswidth} 
\newlength{\DIFdelgraphicsheight} 
\LetLtxMacro{\DIFOincludegraphics}{\includegraphics} 
\newcommand{\DIFaddincludegraphics}[2][]{{\color{blue}\fbox{\DIFOincludegraphics[#1]{#2}}}} 
\newcommand{\DIFdelincludegraphics}[2][]{
\sbox{\DIFdelgraphicsbox}{\DIFOincludegraphics[#1]{#2}}
\settoboxwidth{\DIFdelgraphicswidth}{\DIFdelgraphicsbox} 
\settoboxtotalheight{\DIFdelgraphicsheight}{\DIFdelgraphicsbox} 
\scalebox{\DIFscaledelfig}{
\parbox[b]{\DIFdelgraphicswidth}{\usebox{\DIFdelgraphicsbox}\\[-\baselineskip] \rule{\DIFdelgraphicswidth}{0em}}\llap{\resizebox{\DIFdelgraphicswidth}{\DIFdelgraphicsheight}{
\setlength{\unitlength}{\DIFdelgraphicswidth}
\begin{picture}(1,1)
\thicklines\linethickness{2pt} 
{\color[rgb]{1,0,0}\put(0,0){\framebox(1,1){}}}
{\color[rgb]{1,0,0}\put(0,0){\line( 1,1){1}}}
{\color[rgb]{1,0,0}\put(0,1){\line(1,-1){1}}}
\end{picture}
}\hspace*{3pt}}} 
} 
\LetLtxMacro{\DIFOaddbegin}{\DIFaddbegin} 
\LetLtxMacro{\DIFOaddend}{\DIFaddend} 
\LetLtxMacro{\DIFOdelbegin}{\DIFdelbegin} 
\LetLtxMacro{\DIFOdelend}{\DIFdelend} 
\DeclareRobustCommand{\DIFaddbegin}{\DIFOaddbegin \let\includegraphics\DIFaddincludegraphics} 
\DeclareRobustCommand{\DIFaddend}{\DIFOaddend \let\includegraphics\DIFOincludegraphics} 
\DeclareRobustCommand{\DIFdelbegin}{\DIFOdelbegin \let\includegraphics\DIFdelincludegraphics} 
\DeclareRobustCommand{\DIFdelend}{\DIFOaddend \let\includegraphics\DIFOincludegraphics} 
\LetLtxMacro{\DIFOaddbeginFL}{\DIFaddbeginFL} 
\LetLtxMacro{\DIFOaddendFL}{\DIFaddendFL} 
\LetLtxMacro{\DIFOdelbeginFL}{\DIFdelbeginFL} 
\LetLtxMacro{\DIFOdelendFL}{\DIFdelendFL} 
\DeclareRobustCommand{\DIFaddbeginFL}{\DIFOaddbeginFL \let\includegraphics\DIFaddincludegraphics} 
\DeclareRobustCommand{\DIFaddendFL}{\DIFOaddendFL \let\includegraphics\DIFOincludegraphics} 
\DeclareRobustCommand{\DIFdelbeginFL}{\DIFOdelbeginFL \let\includegraphics\DIFdelincludegraphics} 
\DeclareRobustCommand{\DIFdelendFL}{\DIFOaddendFL \let\includegraphics\DIFOincludegraphics} 

\begin{document}
\title{\Large\bf Oriented Diameter of Mixed Graphs with Given Maximum Undirected Degree}
\author{Ran An$^{a}$, Hengzhe Li$^{a}$, Jianbing Liu\thanks{Corresponding author}$^{*,b}$, Gaoxing Sun$^{a}$\\
\small $^{a}$School of Mathematics and Statistics,\\
\small Henan Normal University, Xinxiang 453007, P.R. China\\
\small $^b$ Department of Mathematics, \\
\small University of Hartford, West Hartford 06117, USA\\
\small Email: rananhtu@163.com, lihengzhe@htu.edu.cn, jianliu@hartford.edu, sungaoxing@htu.edu.cn}
\date{}
\maketitle
\begin{abstract}
In 2018, Dankelmann, Gao, and Surmacs [J. Graph Theory, 88(1): 5--17, 2018] established sharp bounds on the oriented diameter of a bridgeless undirected graph and a bridgeless undirected bipartite graph in terms of vertex degree. In this paper, we extend these results to \emph{mixed graphs}, which contain both directed and undirected edges.

Let the \emph{undirected degree} $d^*_G(x)$ of a vertex $x \in V(G)$ be the number of its incident undirected edges in a mixed graph $G$ of order $n$, and let the \emph{maximum undirected degree} be $\Delta^*(G) = \max\{d^*_G(v) : v \in V(G)\}$. We prove that
\begin{align*}
\text{(1)}\quad & \overrightarrow{\mathrm{diam}}(G) \leq n - \Delta^* + 3 && \text{if $G$ is undirected, or contains a vertex $u$ with $d^*_G(u) = \Delta^*$} \\
& && \text{and $d^+_G(u) + d^-_G(u) \geq 2$, or $\Delta^* = 5$ and $d^+_G(u) + d^-_G(u) = 1$;} \\
\text{(2)}\quad & \overrightarrow{\mathrm{diam}}(G) \leq n - \Delta^* + 4 && \text{otherwise}.
\end{align*}
We also establish bounds for mixed bipartite graphs. If $G$ is a bridgeless mixed bipartite graph with partite sets $A$ and $B$, and $u \in B$, then
\begin{align*}
\text{(1)}\quad & \overrightarrow{\mathrm{diam}}(G) \leq 2(|A| - d(u)) + 7 && \text{if $G$ is undirected;\ \ \ \ \ \ \ \ \ \ \ \ \ \ \ \ \ \ \ \ \ \ \ \ \ \ \ \ \ \ \ \ \ \ \ \ \ \ \ \ \ \ } \\
\text{(2)}\quad & \overrightarrow{\mathrm{diam}}(G) \leq 2(|A| - d^*(u)) + 8 && \text{if $d^+_G(u) + d^-_G(u) \geq 2$;} \\
\text{(3)}\quad & \overrightarrow{\mathrm{diam}}(G) \leq 2(|A| - d^*(u)) + 10 && \text{otherwise}.
\end{align*}

All of the above bounds are sharp, except possibly the last one.

\end{abstract}

\section{Introduction}
A \emph{mixed graph} $G$ is an ordered pair $(V(G), E(G))$, where $V(G)$ is the vertex set and $E(G)$ includes both undirected and directed edges. In this paper, we consider only mixed graphs whose underlying graph is simple.

An \emph{orientable path} from $x_1$ to $x_k$ in a mixed graph $G$ is a sequence of distinct vertices $x_1, x_2, \ldots, x_k$ such that for each $1 \le i \le k-1$, either $x_ix_{i+1} \in E(G)$ (undirected) or $\overrightarrow{x_ix_{i+1}} \in E(G)$ (directed). We denote such a path by $P_{x_1x_k}$.
 An \emph{orientable cycle} is a closed orientable path. For an undirected graph $G$, the orientable paths and cycles are simply paths and cycles. For a directed graph $G$, the orientable paths and cycles are simply directed paths and cycles. A path with $k$ vertices is called a \emph{$k$-path}, while a cycle with $k$ vertices is called a \emph{$k$-cycle}.

A mixed graph $G$ is \emph{connected} if, for any two vertices $x, y \in V(G)$, there exists an orientable path from $x$ to $y$ and an orientable path from $y$ to $x$. In particular, a directed graph is \emph{strongly connected} if it is connected and every edge is a directed edge. A \emph{bridge} in a mixed graph $G$ is an undirected edge $e$ such that $G - e$ is disconnected. A mixed graph is \emph{bridgeless} if it is connected and has no bridge.

The \emph{distance} from $x$ to $y$ in a mixed graph $G$ is the number of edges in a shortest orientable path from $x$ to $y$, and is denoted by $d_G(x, y)$. The \emph{diameter} of $G$ is $\max\{d_G(x, y) : x, y \in V(G)\}$, denoted by $\mathrm{diam}(G)$.

For an orientable path $Q = x_1 \ldots x_k$ and any two integers $1 \le i \le j \le k$, we use $Q(x_i, x_j)$ to denote the orientable subpath $x_i \ldots x_j$. In a mixed graph $G$, for two disjoint paths $Q = x_1 \ldots x_k$ and $R = y_1 \ldots y_\ell$ such that $x_k y_1 \in E(G)$ is either an undirected edge or a directed edge from $x_k$ to $y_1$, we use $QR$ to denote the path $x_1 \ldots x_k y_1 \ldots y_\ell$. We simply write $x_1 R$ rather than $QR$ if $Q = x_1$, and write $Q y_1$ rather than $QR$ if $R = y_1$.

In a mixed graph $G$, for each vertex $x$, define the following neighbor sets:
\[
 N_G^+(x) = \{ y : \overrightarrow{xy} \in E(G) \},
 N_G^-(x) = \{ y : \overrightarrow{yx} \in E(G) \},
 N_G^*(x) = \{ y : xy \in E(G) \},
\]
\[
 N_G(x) = N_G^+(x) \cup N_G^-(x) \cup N_G^*(x),
 N_G[x] = N_G(x) \cup \{ x \}.
\]
The \emph{degree} of $x$ (i.e., the out-degree, in-degree, and undirected degree, respectively) are defined as follows:
\[
d_G(x) = |N_G(x)|, \quad d_G^+(x) = |N_G^+(x)|, \quad d_G^-(x) = |N_G^-(x)|, \quad d_G^*(x) = |N_G^*(x)|.
\]
The \emph{maximum undirected degree} $\Delta^*(G)$ is $\max\{ d_G^*(v) : v \in V(G) \}$. An \emph{$(n, \Delta^*)$-graph} is a mixed graph with $n$ vertices and maximum undirected degree $\Delta^*$.

Assigning a direction to all undirected edges in a mixed graph $G$ is called an \emph{orientation} of $G$. An orientation is \emph{strong} if it is strongly connected. A mixed graph $G$ is \emph{strongly orientable} if it has a strong orientation.
The orientations of two subgraphs of a mixed graph are \emph{compatible} if all their common edges are oriented in the same direction.
The \emph{oriented diameter} $\overrightarrow{\mathrm{diam}}(G)$ is defined as
\[
\overrightarrow{\mathrm{diam}}(G) = \min\{ \mathrm{diam}(\overrightarrow{G}) : \overrightarrow{G} \text{ is a strong orientation of } G \}.
\]

Robbins's Theorem~\cite{Robbins} stated that a graph $G$ has a strong orientation if and only if $G$ is a bridgeless connected graph. According to this theorem, many researchers have studied the oriented diameter of bridgeless graphs. In 1978, Chvátal and Thomassen~\cite{Chvatal} proved that there exists a function $f$ such that, for any undirected bridgeless graph with diameter $d$, its oriented diameter is at most $f(d)$. Specifically, they showed that $\frac{1}{2}d^2 + d \le f(d) \le 2d^2 + 2d$. In 2010, Kwok, Liu, and West~\cite{Kwok} improved this range by proving that $9 \le f(3) \le 11$. Finally, in 2022, Wang and Chen~\cite{Wang2} established that $f(3) = 9$.

Bridgeless mixed graphs can be studied similarly to bridgeless undirected graphs. In 1985, Chung, Garey, and Tarjan~\cite{Chung} proved that a mixed graph $G$ has a strong orientation if and only if $G$ is a bridgeless connected graph, and also showed that for a bridgeless mixed graph $G$ with diameter $d$, its oriented diameter is at most $8d^2 + 8d$.
In 2023, Babu, Benson, and Rajendraprasad~\cite{Babu} further showed that there exists a function $f_m$ such that any bridgeless mixed graph with diameter $d$ has oriented diameter at most $f_m(d)$. Specifically, they proved that $f_m(d) \le 3d^2 + 2d + 2.$

Bounds on $\overrightarrow{\mathrm{diam}}(G)$ have been studied with respect to various graph parameters, including minimum degree~\cite{Bau, Cochran, Czabarka, Huang, Surmacs}, maximum degree~\cite{Dankelmann, Chen}, domination number~\cite{Fomin, Kurz, Dankelmann}, and specific graph classes~\cite{Huang, Kumar, Mondal, Wang1}.

For an undirected graph $G$ of order $n$ and maximum degree $\Delta$, it is known (e.g.,~\cite{Bosak}) that $\mathrm{diam}(G) \le n - \Delta + 1$. In 2018, Dankelmann, Gao, and Surmacs~\cite{Babu} proved that

\begin{theorem}[\upshape Dankelmann, Gao and Surmacs \cite{Dankelmann}]\label{Dmaximumdegree}
For each bridgeless undirected graph $G$ of order $n$ and maximum degree $\Delta$,
$\overrightarrow{diam}(G)\leq n-\Delta+3,$
\end{theorem}

\begin{theorem}[\upshape Dankelmann, Gao and Surmacs \cite{Dankelmann}]\label{Dbipartite}
Let $G=(V(G), E(G))$ be a bridgeless bipartite graph with partite sets $A$ and $B$. For each vertex $u\in B$,
$\overrightarrow{diam}(G)\le 2(|A|-deg_G(u))+7.$
\end{theorem}

This paper extends the results of Dankelmann, Gao and Surmacs to bridgeless mixed graphs through the following findings.

\begin{theorem}\label{maximumdegree}
Let $G$ be a bridgeless mixed graph of order $n$. For any $u \in V(G)$,
\[
\overrightarrow{\mathrm{diam}}(G) \le
\begin{cases}
n - d_G^*(u) + 3, & \text{if $G$ is undirected}, \\
                 & \text{or if } d_G^+(u) + d_G^-(u) \ge 2, \\
                 & \text{or if } d_G^*(u) = 5 \text{ and } d_G^+(u) + d_G^-(u) = 1;\ \ \ \ \ \ \ \ \ \ \ \ \ \ \ \ \ \ \ \ \ \ \ \ \ \ \ \  \ \ \  \\
n - d_G^*(u) + 4, & \text{otherwise.}
\end{cases}
\]
\end{theorem}

\begin{corollary}\label{cor1}
For any bridgeless $(n, \Delta^*)$-graph $G$,
\[
\overrightarrow{\mathrm{diam}}(G) \le
\begin{cases}
n - \Delta^* + 3, & \text{if $G$ is undirected}, \\
                 & \text{or if $G$ has a vertex $u$ with $d_G^*(u) = \Delta^*$ and } d_G^+(u) + d_G^-(u) \ge 2, \\
                 & \text{or if $G$ has a vertex $u$ with $d_G^*(u) = \Delta^* = 5$ and } d_G^+(u) + d_G^-(u) = 1; \\
n - \Delta^* + 4, & \text{otherwise}.
\end{cases}
\]
\end{corollary}

\begin{theorem}\label{bipartite}
Let $G$ be a bridgeless mixed bipartite graph of order $n$ with partite sets $A$ and $B$. For any $u \in B$,
\[
\overrightarrow{\mathrm{diam}}(G) \le
\begin{cases}
2(|A| - d(u)) + 7, & \text{if $G$ is undirected; \ \ \ \ \ \ \ \ \ \ \ \ \ \ \ \ \  \ \ \ \ \ \ \ \ \ \ \ \ \  \ \ \  \ \ \ \ \ \ \ \ \ \ \ \ \ \ \ \  } \\
2(|A| - d_G^*(u)) + 8, & \text{if } d_G^+(u) + d_G^-(u) \ge 2; \\
2(|A| - d_G^*(u)) + 10, & \text{otherwise}.
\end{cases}
\]
\end{theorem}

\begin{corollary}\label{equal}
For any bridgeless mixed equal-partition bipartite graph $G$ of order $n$ with maximum undirected degree $\Delta^*$,
\[
\overrightarrow{\mathrm{diam}}(G) \le
\begin{cases}
n - 2\Delta^* + 7, & \text{if $G$ is undirected}; \\
n - 2\Delta^* + 8, & \text{if $G$ has a vertex $u$ with $d_G^*(u) = \Delta^*(G)$ and } d_G^+(u) + d_G^-(u) \ge 2; \\
n - 2\Delta^* + 10, & \text{otherwise}.
\end{cases}
\]
\end{corollary}

This paper is structured as follows: Section~2 introduces key lemmas supporting the main results. Section~3.1 provides proofs for Theorem~1.3, along with sharp examples. Section~3.2 provides proofs for Theorem~1.5, along with sharp examples, except when $G$ is a mixed equal-partition bipartite graph with at least one directed edge and $d_G(u) \le d_G^*(u) + 1$ for each vertex $u$ with maximum undirected degree.

\section{Key Lemmas}

This section introduces key lemmas that support our main results. Let $G$ be a bridgeless mixed graph. Fix a vertex $u \in V(G)$. We construct a mixed graph $G_1$ from $G$ by orienting certain undirected edges incident to $u$. Next, we form a subgraph $H_u \subseteq G_1$ such that $N_G(u) \subseteq V(H_u)$, and we orient $H_u$ into a strongly connected graph $\overrightarrow{H_u}$. Finally, within $\overrightarrow{H_u}$, we construct directed paths $T_{uv}$ from $u$ to $v$ and $T_{vu}$ from $v$ to $u$, each containing as few vertices from $N_G(u)$ as possible.
These objectives are achieved in stages.

\textbf{Stage 1.} In this stage, we divide $N_G(u)$ into three parts and orient some of the undirected edges of $G$ incident with $u$.

For a vertex $u$ of a bridgeless mixed graph $G$, we divide $N_G(u)$ into three parts:
\[
X_1 = N_G^-(u), \quad Y_1 = N_G^+(u), \quad \text{and} \quad N_{\mathrm{un}} = N_G^*(u).
\]
For every $v \in N_G(u)$, let $\ell_G(v)$ denote the length of a shortest orientable cycle in $G$ containing the edge between $u$ and $v$. Define
\[
s = \sum_{v \in N_G(u)} \ell_G(v).
\]
We complete the orientation process in several stages. We now orient some of the undirected edges of $G$ incident with $u$. We initialize $X_2 = Y_2 = Z = \emptyset$.

Repeat steps (1)--(3) for all vertices in $N_{\mathrm{un}}$:

\begin{enumerate}
    \item [(1)] The edge $uv$ is oriented from $v$ to $u$ if the parameter $s$ remains unchanged after this orientation. In this case, the vertex $v$ is added to $X_2$.

    \item [(2)] Otherwise, the edge $uv$ is oriented from $u$ to $v$ if the parameter $s$ remains unchanged after this orientation. The vertex $v$ is then added to $Y_2$.

    \item [(3)] Otherwise, the edge $uv$ remains undirected. Such an edge $uv$ is called a \emph{conflicted edge}, and the vertex $v$ is called a \emph{conflicted vertex}. In this case, $v$ is added to $Z$.
\end{enumerate}

Let $X = X_1 \cup X_2$ and $Y = Y_1 \cup Y_2$.

Babu, Benson, and Rajendraprasad~\cite{Bau} studied unavoidable edges of the shortest orientable cycle containing conflict edges. We use an improved version of their result, due to Li, Ding, Liu et al.~\cite{Li}.

\begin{lemma}[Li, Ding, Liu et al.~\cite{Li}]\label{Fix}
If an edge $uv_i$ is conflicted, then one of the following must hold:
\begin{itemize}
    \item [\rm (1)] There exists a vertex $v_j \in X$ such that $uv_i$ is included in every shortest orientable cycle containing $\overrightarrow{v_j u}$.

    \item [\rm (2)] There exists a vertex $v_k \in Y$ such that $uv_i$ is included in every shortest orientable cycle containing $\overrightarrow{u v_k}$.
\end{itemize}
\end{lemma}

Let $G_1$ be the mixed graph derived from $G$ by following the orientation described above. For $G_1$, the following observations hold.

\begin{observation}\label{shortestcycle}
For each $v \in N_G(u)$, $\ell_{G_1}(v) = \ell_G(v)$.
\end{observation}

\begin{observation}\label{internalvertex}
The following properties hold for shortest orientable paths between $u$ and its neighbors in $G_1$:
\begin{enumerate}
    \item [\rm (1)] For each $x \in X$, the second vertex of each shortest orientable path $P_{ux}$ in $G_1$ from $u$ to $x$ in $G$ belongs to $Y \cup Z$.

    \item [\rm (2)] For each $y \in Y$, the penultimate vertex of each shortest orientable path $P_{yu}$ from $y$ to $u$ in $G_1$ belongs to $X \cup Z$.
\end{enumerate}
\end{observation}

\begin{lemma}\label{path}
If  $G_1$ is constructed as described above, then the following properties hold.

\begin{enumerate}
    \item [\rm (1)] For each $x \in X$, each shortest orientable path $P_{ux}$ from $u$ to $x$ in $G_1$ contains exactly one vertex from $Y \cup Z$ and no vertex from $X_2 \setminus \{x\}$.

    \item [\rm (2)]  For each $y \in Y$, each shortest orientable path $P_{yu}$ from $y$ to $u$ in $G_1$ contains exactly one vertex from $X \cup Z$ and no vertex from $Y_2 \setminus \{y\}$.
\end{enumerate}
\end{lemma}

\begin{proof}
Let $x \in X$, and let $P_{ux}$ be a shortest orientable path in $G_1$ from $u$ to $x$. Suppose that $P_{ux}$ contains at least two vertices from $Y \cup Z$. Then, we can select two vertices $u_1, u_2 \in V(P_{ux}) \cap (Y \cup Z)$ such that the subpath $P_{ux}(u, u_2)$ contains $u_1$. In this case, the path $uP_{ux}(u_2, x)$ forms an orientable path in $G_1$ that is shorter than $P_{ux}$, leading to a contradiction.

Alternatively, if $P_{ux}$ contains a vertex, say $u_1$, from $X_2$ other than $x$, then the path $uP_{ux}(u_1, x)u$ forms an orientable cycle in $G$ that is shorter than the cycle $P_{ux}u$ in $G_1$, again contradicting Observation 2.4. This proves part (1).

By the symmetry of $X$ and $Y$, the same argument establishes part (2).
\end{proof}

\textbf{Stage 2.} At this stage, we construct a strongly connected subgraph $\overrightarrow{H_u}$ that covers all vertices in $N_G(u)$, along with directed paths $T_{uv}$ from $u$ to $v$ and $T_{vu}$ from $v$ to $u$, each containing as few vertices from $N_G(u)$ as possible.

By the symmetry of $X$ and $Y$, we can obtain an in-tree $T_{\mathrm{in}}$ rooted at $u$ with leaf set $Y$.

The orientable paths $P_{ux}$ and $P_{yu}$ constructed by Algorithm~1 satisfy the following lemma.

\begin{algorithm}[H]
\caption{Out-tree $T_{\text{out}}$ Construction}
\textbf{Input:} A strongly connected bridgeless mixed graph $G$ and a vertex $u \in V(G)$.\\
\textbf{Output:} $T_{\text{out}}$ is an out-tree rooted at $u$ with leaf set $X$ and a shortest orientable path $P_{ux}$ from $u$ to $x$ for each $x \in X$ in $T_{\text{out}}$.
\begin{algorithmic}[1]
    \State Perform a breadth-first search in $G_1$ starting from vertex $u$ to construct a spanning tree $T$.
    \State Initialize $T_{\text{out}} \gets \{u\}$.
    \While{there exists a leaf $x$ in $T$ such that $x \in X$}
        \State Add $x$ and the oriented path from $u$ to $x$ in $T$ to $T_{\text{out}}$.
        \State Let $P_{ux}$ be the orientable path from $u$ to $x$ in $T_{\text{out}}$.
    \EndWhile
\end{algorithmic}
\end{algorithm}

\begin{lemma}\label{vertexnumber1}
Assume that the paths $P_{ux}$ and $P_{yu}$ are constructed as in Stage~2.
\begin{enumerate}
    \item [\rm (1)] For each $x \in X$, the path $P_{ux}$ contains exactly one vertex from $Y \cup Z$ and no vertex from $X_2 \setminus \{x\}$.

    \item [\rm (2)] For each $y \in Y$, the path $P_{yu}$ contains exactly one vertex from $X \cup Z$ and no vertex from $Y_2 \setminus \{y\}$.
\end{enumerate}
\end{lemma}
\begin{proof} From Steps~1 to~5 of Algorithm~1, we conclude that $P_{ux}$ is a shortest orientable path from $u$ to $x$ in both $G$ and $G_1$ for each $x \in X$. By Lemma~2.4, (1) follows.
By the symmetry of $X$ and $Y$, (2) can be proved similarly.
\end{proof}

From Lemma \ref{vertexnumber1}, we obtain the following result.

\begin{lemma}\label{vertexnumber1+}
For each $x \in X$, the path $P_{ux}$ contains at most $d_G^-(u) + 2$ vertices from $N_G(u)$.
\end{lemma}

By the process of constructing $P_{ux}$ and $P_{yu}$, the following observations hold.

\begin{observation}\label{shortestpath}
The constructed paths $P_{ux}$ and $P_{yu}$ have the following properties.
\begin{enumerate}
    \item [\rm (1)] For any $x, x' \in X$, if $P_{ux}$ and $P_{ux'}$ share a vertex $v$, then the subpaths $P_{ux}(u, v)$ and $P_{ux'}(u, v)$ are identical.

    \item [\rm (2)] For each $y, y' \in Y$, if $P_{yu}$ and $P_{y'u}$ have a common vertex $w$, then the subpaths $P_{yu}(w, u)$ and $P_{y'u}(w, u)$ are identical.
\end{enumerate}
\end{observation}

Assume $X = \{x_1, x_2, \ldots, x_s\}$ and $Y = \{y_1, y_2, \ldots, y_t\}$. Based on Algorithm~1, for each $x_i \in X$, the path $P_{ux_i}$ is a shortest orientable path in both $T_{\text{out}}$ and $G_1$ from $u$ to $x_i$; for each $y_j \in Y$, the path $P_{y_ju}$ is a shortest orientable path in both $T_{\text{in}}$ and $G_1$ from $y_j$ to $u$.
Define $C_{x_i} := x_i P_{ux_i}$, $C_{y_j} := P_{y_j} y_j$, and let $H_X := \bigcup_{i=1}^{s} C_{x_i}$. For each $C_{x_i}$, orient it to form the directed cycle $\overrightarrow{C_{x_i}}$, and set $\overrightarrow{H_X} := \bigcup_{i=1}^{s} \overrightarrow{C_{x_i}}$.
By Observation \ref{Fix}, we know that
\[
N_{G_1}(u) \subseteq \bigcup_{i=1}^{s} V(C_{x_i}) \cup \bigcup_{j=1}^{t} V(C_{y_j}).
\]

Using $H_X$ and $C_{y_j}$, Algorithm~2 constructs a mixed subgraph $H_u$ of $G_1$ such that $N_G(u) \subseteq V(H_u)$, its strong orientation $\overrightarrow{H_u}$, and an orientable path $Q_{yu}$ for every $y \in Y$.

\begin{algorithm}[htb]
\renewcommand{\baselinestretch}{1.3}
\caption{Strong Orientation and Path Construction}
\begin{algorithmic}[1]  
\Require A strongly connected bridgeless mixed graph $H_X$ and the orientable path $P_{yu}$ for each $y\in Y$.
\Ensure A strongly connected bridgeless graph $\overrightarrow{H_u}$ with $N_G[u]\in V(H_u)$ and the directed path $Q_{yu}$ from $y$ to $u$ for each $y\in Y$ in $\overrightarrow{H_u}$.
\State Initialize $j:=0$, $Q_{y_ju}:=P_{y_ju}$, $H_u:=H_X=\bigcup^s_{i=1}C_{x_i}$, and $\overrightarrow{H_u}:=\overrightarrow{H_X}= \bigcup_{i=1}^s \overrightarrow{C_{x_i}}$

\While{$j < t$}
    \State $j:=j+1$.
    \If{$y_j \in Y \cap V(H_X)$}
        \State Pick $x \in X$ such that $P_{ux}$ contains $y_j$ and $V(C_x) \cap X = \{x\}$.
        \State Let $Q_{y_ju}$ be the mixed path in $C_x$ from $y_j$ to $u$.
    \Else\ {$y_j \in Y \setminus V(H_X)$}
        \If {$V(P_{y_ju}) \cap V(H_X) = \emptyset$}
            \State $Q_{y_ju} := P_{y_ju}$, $C^1_{y_j} := Q_{y_ju}y_j$, and $C^2_{y_j} = \emptyset$.
             \State  $H_u:=H_u\cup C^1_{y_j}$ and $\overrightarrow{H_u}:=\overrightarrow{H_u}\cup \overrightarrow{C^1_{y_j}}$.
        \Else\ {$V(P_{y_ju}) \cap V(H_X) \neq \emptyset$}
            \State  Let $u_{y_j}$ and $v_{y_j}$ be the first and last vertices of $P_{y_ju}$ in $H_X$, respectively.
            \State Choose $x \in X$ such that $P_{ux}$ contains $u_{y_j}$ and $P_{ux}(u_{y_j},x)$ is minimal.
            \State Choose $x'\in X$ such that $v_{y_j}$ is contained in $C_{x'}$.

            \If{$u_{y_j} = v_{y_j} \in Z$}
                \State $Q_{y_ju} := P_{y_ju}(y_j, u_{y_j})P_{ux}(u_{y_j}, x)u$, $C^1_{y_j} := Q_{y_ju}y_j$, and $C^2_{y_j} := \emptyset$.
                \State $H_u:=H_u\cup C^1_{y_j}$ and $\overrightarrow{H_u}:=\overrightarrow{H_u}\cup \overrightarrow{C^1_{y_j}}$.
            \ElsIf {$u_{y_j} = v_{y_j} \notin Z$}
                \State $Q_{y_ju} := P_{y_ju}$, $C^1_{y_j} := Q_{y_ju}y_j$, and $C^2_{y_j} := \emptyset$.
                \State $H_u:=H_u\cup C^1_{y_j}$ and $\overrightarrow{H_u}:=\overrightarrow{H_u}\cup \overrightarrow{C^1_{y_j}}$.
            \Else\ {$u_{y_j} \neq v_{y_j}$}
                \State $Q_{y_ju} := P_{y_ju}(y_j, u_{y_j})P_{ux}(u_{y_j}, x)u$, $C^1_{y_j} := Q_{y_ju}y_j$, and
                 $C^2_{y_j} := P_{ux'}(u, v_{y_j})P_{y_ju}(v_{y_j},u)$.
                \State $H_u:=H_u\cup C^1_{y_j}\cup C^2_{y_j}$ and $\overrightarrow{H_u}:=\overrightarrow{H_u}\cup \overrightarrow{C^1_{y_j}}\cup \overrightarrow{C^2_{y_j}}$.
            \EndIf
        \EndIf
    \EndIf
\EndWhile
\end{algorithmic}
\end{algorithm}

One has the following lemma.
\begin{lemma}\label{vertexnumber2}
The directed path $Q_{yu}$ constructed by Algorithm~2 satisfies the following properties.
\begin{enumerate}
    \item [\rm (1)] For each $y \in Y$, $Q_{yu}$ contains no vertex from $Y_2 \setminus \{y\}$.
    If $Z = \emptyset$, then $Q_{yu}$ contains exactly one vertex from $X$.
    If $Z \ne \emptyset$, then $Q_{yu}$ contains at most two vertices from $X \cup Z$.
    \item [\rm (2)]  For each $y \in Y$ such that $C^2_y \ne \emptyset$, the orientable subpath of $C^2_y$ from $u$ to $Z$ contains no vertices from $X_2$.
\end{enumerate}
\end{lemma}

\begin{proof}
First, we show (1). By Lemma \ref{vertexnumber1}, it follows that for each $y \in Y$, the path $P_{yu}$ contains exactly one vertex from $X \cup Z$ and no vertex from $Y_2 \setminus \{y\}$. If $Q_{yu} := P_{yu}$, the result holds directly. Otherwise, for a vertex $y \in Y$, let $Q_{yu} := P_{yu}(y, u_y) P_{ux}(u_y, x) u$, where $x \in X$ is such that $u_y \in V(C_x)$ and the orientable path from $u_y$ to $u$ in $C_x$ contains exactly one vertex from $X$. Since $P_{yu}(y, u_y)$ is a subpath of $P_{yu}$, it contains no vertices from $Y_2$. Moreover, since $P_{ux}(u_y, x)$ is a subpath of $P_{ux}$ and $u_y \in P_{yu}$, by Lemma~2.5, $u_y \in Y_1$. Therefore, $Q_{yu} := P_{yu}(y, u_y) P_{ux}(u_y, x) u$ contains no vertex from $Y_2 \setminus \{y\}$.  By the definition of $Q_{yu}$, it is also clear that $Q_{yu}$ contains at most two vertices from $X \cup Z$ if $Z \ne \emptyset$, and exactly one vertex from $X$ if $Z = \emptyset$.

Now, we show (2). For any $y \in Y$, let $v_y$ be the last vertex of $P_{yu}$ that lies in $C'_x$, and define $C^2_y := P_{ux'}(u, v_y) P_{yu}(v_y, u)$. Since $P_{yu}(v_y, u)$ is a subpath of $P_{yu}$, Lemma~2.5 implies that it contains no vertices from $X_2$. If $P_{ux'}(u, v_y)$ contains vertices from $X_2$, then $v_y \in X_2$, contradicting the assumption. In that case, we could construct a shorter orientable cycle $u P_{yu}(y, v_y) u$ than $P_{yu}$, which contradicts minimality. Hence, $C^2_y$ contains no vertex from $X_2$, and the claim is proved.
\end{proof}

From Lemma \ref{vertexnumber2}, we obtain the following result.

\begin{lemma}\label{vertexnumber2+}
For each $y \in Y$, the path $Q_{yu}$ contains at most $d_G^+(u) + 3$ vertices from $N_G(u)$ in $\overrightarrow{H_u}$. Furthermore, if $Z = \emptyset$, then $Q_{yu}$ contains at most $d_G^-(u) + 2$ vertices from $N_G(u)$ in $\overrightarrow{H_u}$.
\end{lemma}

Next, we use Algorithm~3 to construct a directed path $Q_{uz}$ from $u$ to $z$ and a directed path $Q_{zu}$ from $z$ to $u$ for each $z \in Z$.

\begin{algorithm}
\renewcommand{\baselinestretch}{1.3}
\caption{{Directed Path Construction for Conflicted Vertices in $\overrightarrow{H_u}$}}
\begin{algorithmic}[1]  
\Require A strongly connected bridgeless graph $\overrightarrow{H_u}$.
\Ensure A directed path $Q_{uz}$ from $u$ to $z$ in $\overrightarrow{H_u}$ and a directed path $Q_{zu}$ from $z$ to $u$ in $\overrightarrow{H_u}$ for each $z\in Z$.

\For{each $z \in Z$}
    \If{there exists a vertex $x\in X$ such that $z \in Z \cap V(\overrightarrow{C_x})$}
        \State Choose the first vertex $x'$ of the directed path $\overrightarrow{C_x}(z,u)$ belonging to $X$.
        \State $Q_{uz} := uz$ and $Q_{zu} := \overrightarrow{C_x}(z,x')u$.
    \ElsIf{there exists a vertex $y\in Y$ such that $z \in Z \cap \bigcup_{y \in Y} V(C^1_y)$}
        \State Choose the last vertex $y'$ of $P_{yu}$ belonging to $Y$.
        \State $Q_{uz} := P_{yu}(y',z)$ and $Q_{zu} := zu$.
    \Else{ Choose a vertex $y\in Y$ such that $z \in Z \cap \bigcup_{y \in Y} V(C^2_y)$}
        \State Let $v_y$ be the last vertex of $P_{yu}$ such that $v_y\in V(H_X)$ and let $x'\in X$ such that $v_y$ lies on $P_{ux'}$.
        \State $\overrightarrow{C^2_y}:= P_{ux'}(u, v_y)P_{yu}(v_y, u)$.

        \If{$P_{yu}(v_y, u)$ includes certain vertices from $Y$}
            \State Pick the last vertex $y'$ of $P_{yu}(v_y, u)$ that belongs to $Y$.
            \State $Q_{uz} := \overrightarrow{C^2_y}(y', z)$ and $Q_{zu} := zu$.

        \Else\  $Q_{uz} := \overrightarrow{C^2_y}(u, z)$ and $Q_{zu} := zu$.
        \EndIf
    \EndIf
\EndFor
\end{algorithmic}
\end{algorithm}

From Step 13 of Algorithm 3, the following observations hold.

\begin{observation}\label{Q_{uz}}
For each $z \in Z$, if $Q_{uz}$ contains vertices from $X_1$, then:
\begin{enumerate}
    \item [\rm (1)] There exists a vertex $y \in Y$ such that $z$ lies on $P_{yu}$.

    \item [\rm (2)] Let $v_y$ be the last vertex of $P_{yu}$ in $H_X$, and let $x' \in X$ be such that $P_{ux'}$ passes through $v_y$. Then
    \[
    Q_{uz} = \overrightarrow{C^2_y}(u, z), \quad \text{where } C^2_y = P_{ux'}(u, v_y) P_{yu}(v_y, u).
    \]
\end{enumerate}
\end{observation}

We have the following lemma.
\begin{lemma}\label{dvertexnumber3}
The directed paths $Q_{uz}$ and $Q_{zu}$ constructed by Algorithm~3 satisfy the following properties.
\begin{enumerate}
    \item [\rm (1)]  For each $z \in Z$, there exists a directed path $Q_{uz}$ from $u$ to $z$ that contains at most $d_G^-(u) + 2$ vertices from $N_G(u)$ in $\overrightarrow{H_u}$.

    \item [\rm (2)]  For each $z \in Z$, there exists a directed path $Q_{zu}$ from $z$ to $u$ that contains at most two vertices from $N_G(u)$ in $\overrightarrow{H_u}$.
\end{enumerate}
\end{lemma}

\begin{proof}
First, we show (1) is true. Based on Algorithm~3, the path $Q_{uz}$ has four possible constructions:
    \begin{itemize}
        \item [\rm (a)] If $Q_{uz} := uz$, then the directed path from $u$ to $z$ contains exactly one vertex from $N_G(u)$.

        \item [\rm (b)] If $Q_{uz} := P_{yu}(y', z)$, where $y'$ is the last vertex of $P_{yu}$ belonging to $Y$, then from Steps~5 to~7 of Algorithm~3, we know that this path contains at most two vertices from $N_G(u)$.

        \item [\rm (c)] If $Q_{uz} := \overrightarrow{C^2_y}(y', z)$, where $y'$ is the last vertex of $P_{yu}(v_y, u)$ that belongs to $Y$, then by Steps~8 to~13 of Algorithm~3, the path from $u$ to $z$ contains at most two vertices from $N_G(u)$.

        \item [\rm (d)] If $Q_{uz} := \overrightarrow{C^2_y}(u, z)$, then from Steps~8 and~14 of Algorithm~3 together with Lemma \ref{vertexnumber2}, the path from $u$ to $z$ contains at most $d_G^-(u) + 2$ vertices from $N_G(u)$.
    \end{itemize}
    Therefore, in all cases, property (1) holds.

Now we shall show (2). Based on Algorithm~3, the path $Q_{zu}$ has two possible cases:
    \begin{itemize}
        \item [\rm (a)] If $Q_{zu} := zu$, then it contains exactly one vertex from $N_G(u)$.
        \item [\rm (b)] Otherwise, we know from Algorithm~3 that $z \in Z \cap V(\overrightarrow{C_x})$ for some $x \in X$, and from Steps~2 to~4 of the algorithm, we have $Q_{zu} := \overrightarrow{C_x}(z, x')u$. This path contains at most two vertices from $N_G(u)$.
    \end{itemize}
Hence, property (2) also holds.
\end{proof}

Let $G_2$ be the mixed graph derived from $G_1$ by following the orientation described above. Since $H_u \subseteq G_1$ consists of a union of mixed cycles, $G_2$ remains a bridgeless strongly connected mixed graph. According to the following theorem, $G_2$ has a strong orientation $\overrightarrow{G}$, which is also a strong orientation of both $G$ and $G_1$.

\begin{theorem}[Chung, Garey, and Tarjan~\cite{Chung}]\label{extend}
A mixed multigraph is strongly orientable if and only if it is strongly connected and has no bridges.
\end{theorem}

To aid the reader in navigating the notation used throughout this section, we summarize below the key symbols and their definitions, see Table 1.
\begin{table}[H]
\centering
\caption{Summary of Symbols Introduced in Section 2}
\renewcommand{\arraystretch}{1.2}
\begin{tabular}{ll}
\hline
\textbf{Symbol} & \textbf{Description} \\
\hline
$N_G(u)$ & The neighborhood of $u$ in graph $G$ \\
$N_G^+(u), N_G^-(u), N_G^*(u)$ & Out-neighbors, in-neighbors, and undirected neighbors of $u$ \\
$X_1 = N_G^-(u)$ & In-neighbors of $u$ \\
$Y_1 = N_G^+(u)$ & Out-neighbors of $u$ \\
$N_{\mathrm{un}} = N_G^*(u)$ & Undirected neighbors of $u$ \\
$X_2, Y_2, Z$ & Partitions of $N_G(u)$ from Stage 1 \\
$X = X_1 \cup X_2$, $Y = Y_1 \cup Y_2$ & Final partition of neighbors \\
$\ell_G(v)$ & Length of the shortest orientable cycle containing $uv$ \\
$H_X$ & Subgraph formed by the union of $C_{x_i}$ cycles for $x_i \in X$ \\
$\overrightarrow{H_X}$ & Strong orientation of $H_X$ \\
$P_{ux}, P_{yu}$ & Orientable paths from $u$ to $x \in X$ and $y \in Y$ to $u$ \\
$T_{\text{out}}, T_{\text{in}}$ & Out-tree and in-tree rooted at $u$ with leaves in $X$ and $Y$ \\
$C_x, C_y$ & Orientable cycles used to define $H_X$ and $H_u$ \\
$C^1_y, C^2_y$ & Auxiliary cycles used to extend $H_X$ to $H_u$ when $y \notin H_X$ \\
$H_u$ & Mixed subgraph of $G_1$ covering $N_G(u)$ \\
$\overrightarrow{H_u}$ & Strong orientation of $H_u$ \\
$Q_{ux}, Q_{yu}$ & Orientable paths in $\overrightarrow{H_u}$ for $x \in X$, $y \in Y$ \\
$Q_{uz}, Q_{zu}$ & Directed paths from $u$ to $z$ and from $z$ to $u$ for $z \in Z$ \\
\hline
\end{tabular}
\end{table}

\section{Proofs of Main Theorems}

This section presents the detailed proofs of our main results stated in Section~1. We begin by proving Theorem \ref{maximumdegree}  and Corollary \ref{cor1}, which establish upper bounds on the oriented diameter of bridgeless mixed graphs. We then prove Theorem \ref{bipartite} and Corollary \ref{equal} concerning the case of bridgeless mixed bipartite graphs. Each proof relies on the constructions and structural lemmas developed in Section~2.

\subsection{The proof of Theorem \ref{maximumdegree}}

In this subsection, we provide a proof of Theorem \ref{maximumdegree} and present a family of examples demonstrating the sharpness of the bound.

\noindent\textbf{The proof of Theorem~\ref{maximumdegree}.}
Let  $u \in V(G) $. By Theorem~\ref{extend}, any strong orientation  $D_u$ of  $H_u$ can be extended to a strong orientation  $D$ of  $G$.

For any two distinct vertices  $v', w' \in V(D)$, choose  $v, w \in V(D_u)$ such that  $d_D(v', v)$ and  $d_D(w, w')$ are minimized. Then
\[
d_D(v', w') \leq d_D(v', v) + d_D(v, w) + d_D(w, w').
\]
To prove Theorem~\ref{maximumdegree}, it suffices to show that  $D_u$ contains a path  $T_{vw}$ from  $v$ to  $w$ that passes through at most  $d^+_G(u) + d^-_G(u) + 4$ vertices in  $N_G(u)$, if any of the following hold:
\begin{itemize}
    \item [(1)]   $G$ is undirected, or
    \item [(2)]   $d^+_G(u) + d^-_G(u) \geq 2$, or
    \item [(3)]   $d^*_G(u) = 5 $ and  $d^+_G(u) + d^-_G(u) = 1$;
\end{itemize}
and at most  $d^+_G(u) + d^-_G(u) + 5$ vertices of  $N_G(u)$, otherwise.

If  $G $ is undirected, then by the argument in the proof of Theorem~\ref{Dmaximumdegree}, $G$ has an our desired strongly orientation. We now consider the following four cases.

\noindent\textbf{Case 1.} $d^+_G(u) = d^-_G(u) = 0$.

In this case, let $D_u := \overrightarrow{H_u}$. Since $X_1 = N^-_G(u) = \emptyset$ and $Y_1 = N^+_G(u) = \emptyset$, the neighborhood of $u$ is partitioned into $X_2$, $Y_2$, and $Z$.

We claim that for each $v \in V(D_u)$, there exists a directed path $T_{vu}$ from $v$ to $u$ containing at most three vertices of $N_G(u)$. Specifically:
\begin{itemize}
  \item [(1)] If $v \in X = X_2$, then $T_{vu} := vu$ contains exactly one vertex of $N_G(u)$;
  \item [(2)] If $v \in Y = Y_2$, then $T_{vu} := Q_{vu}$ contains at most three vertices of $N_G(u)$ by Lemma~\ref{vertexnumber2+};
  \item [(3)] If $v \in Z$, then $T_{vu} := Q_{vu}$ contains at most two vertices of $N_G(u)$ by Lemma~\ref{dvertexnumber3};
  \item [(4)] If $v \in V(D_u) \setminus N_G[u]$, then there exists a vertex $v^\star \in N_G(u)$ such that the subpath of $T_{v^\star u}$ from $v$ to $u$ is a desired path.
\end{itemize}

Similarly, for each $w \in V(D_u)$, there exists a directed path $T_{uw}$ from $u$ to $w$ containing at most two vertices of $N_G(u)$:
\begin{itemize}
  \item [(1)] If $w \in X = X_2$, then $T_{uw} := P_{uw}$ contains at most two vertices of $N_G(u)$ by Lemma~\ref{vertexnumber1+};
  \item [(2)] If $w \in Y$, then $T_{uw} := uw$ contains exactly one vertex of $N_G(u)$;
  \item [(3)] If $w \in Z$, then $T_{uw} := Q_{uw}$ contains at most two vertices of $N_G(u)$ by Lemma~\ref{dvertexnumber3};
  \item [(4)] If $w \in V(D_u) \setminus N_G[u]$, then there exists a vertex $w^\star \in N_G(u)$ such that the subpath of $T_{uw^\star}$ from $u$ to $w$ is a desired path.
\end{itemize}

Therefore, for every pair $v, w \in V(D_u)$, there exists a directed $vw$-path $T_{vw}$ in $D_u$ containing at most five vertices of $N_G(u)$.

\noindent\textbf{Case 2.} $d^+_G(u) + d^-_G(u) = 1$.

By symmetry, we may assume that $d^+_G(u) = 1$ and $d^-_G(u) = 0$. Then the neighborhood of $u$ is partitioned into $X_2$, $Y = Y_1 \cup Y_2$, and $Z$.

We claim that for each $v \in \overrightarrow{H_u}$, there exists a path $T_{vu}$ containing at most four vertices of $N_G(u)$:
\begin{itemize}
    \item[(1)] If $v \in X$, then $T'_{vu} := vu$  contains exactly one vertex of $N_G(u)$;
    \item[(2)] If $v \in Y$, then $T'_{vu} := Q_{vu}$  contains at most four vertices of $N_G(u)$ by Lemma~\ref{vertexnumber2+};
    \item[(3)] If $v \in Z$, then $T'_{vu} := Q_{vu}$  contains at most two vertices of $N_G(u)$ by Lemma~\ref{dvertexnumber3};
    \item[(4)] If $v \in V(D_u) \setminus N_G[u]$, then there exists a vertex $v^\star \in N_G(u)$ such that the subpath of $T'_{v^\star u}$  from $v$ to $u$ is a desired path.
\end{itemize}

We also claim that for each $w \in \overrightarrow{H_u}$, there exists a path $T_{uw}$ containing at most two vertices of $N_G(u)$:
\begin{itemize}
    \item[(1)] If $w \in X$, then $T'_{uw} := P_{uw}$ contains at most two vertices of $N_G(u)$ by Lemma~\ref{vertexnumber1+};
    \item[(2)] If $w \in Y$, then $T'_{uw} := uw$  contains exactly one vertex of $N_G(u)$;
    \item[(3)] If $w \in Z$, then $T'_{uw} := Q_{uw}$  contains at most two vertices of $N_G(u)$ by Lemma~\ref{dvertexnumber3};
    \item[(4)] If $w \in V(\overrightarrow{H_u}) \setminus N_G[u]$, then there exists a vertex $w^\star \in N_G(u)$ such that the subpath of $T'_{uw^\star}$ from $u$ to $w$ is a desired path.
\end{itemize}

If $d^*_G(u) \ge 6$, then let $D_u:=\overrightarrow{H_u}$. 
Therefore, for every $v, w \in V(D_u)$, there exists a $vw$-path $T'_{vw}$ in $D_u$ containing at most six vertices of $N_G(u)$ and we are done.

Suppose that $d^*_G(u) \le 5$. We will reorient an edge (which is an undirected edge in $G$) in $\overrightarrow{H_u}$  so that, for any $v, w \in V(H_u)$, there is a $vw$-path $T_{vw}$ in the resulting graph containing at most five vertices from $N_G(u)$.

If for each $v \in V(\overrightarrow{H_u})$, we have $|T'_{vu} \cap N_G(u)| \le 3$, then $D_u = \overrightarrow{H_u}$ is our desired orientation, and we are done.

Suppose there exists a vertex $s \in V(\overrightarrow{H_u})$ such that $|T'_{su} \cap N_G(u)| = 4$. Furthermore, $T'_{su}$ and $H_X$ contain at least three vertices of $d^*_G(u)$. From Steps 11–16 of Algorithm~2, we conclude that there is exactly one vertex, say $z$, that lies between $T'_{su}$ and $H_X$, and $z \in Z$.

We next show that all the paths $Q_{yu}$ that contain four vertices of $N_G(u)$ must pass through a common conflicted vertex $z \in Z$.

\begin{claim}\label{z}
For each $y \in Y$, if $Q_{yu}$ does not pass through $z$, then it must not pass through any vertices in $Y_1$.
\end{claim}

\begin{proof}[Proof of Claim 1]
We prove this by showing the contrapositive. Suppose there exists a vertex $s' \in Y$ such that $T_{s'u}$ passes through a vertex $y \in Y_1$. By Observation~\ref{shortestpath}, the subpath $T_{s'u}(y, u)$ is identical to the subpath $T_{yu}(y, u)$, and hence $T_{s'u}$ must pass through $z$. Therefore, the contrapositive holds.
\end{proof}

Let $D_u$ be the directed graph obtained from $\overrightarrow{H_u}$ by reversing the arc $\overrightarrow{uz}$ into $\overrightarrow{zu}$.

We claim that for each $v \in V(D_u)$, the path $T_{vu}$ may contain three vertices of $N_G(u)$ if and only if it passes through $z$; otherwise, $T_{vu}$ contains at most two vertices of $N_G(u)$.

\begin{itemize}
    \item[(1)] For $v \in Y$ such that $Q_{vu}$ contains $z$, we construct a path from $v$ to $u$ by replacing the subpath $Q_{vu}(z, u)$ with the arc $zu$. This modified path contains at most three vertices of $N_G(u)$.

    \item[(2)] For $v \in Z$ such that $Q_{vu}$ contains $z$, we let $T_{vu} := vu$, which contains exactly one vertex of $N_G(u)$.

    \item[(3)] Next, we consider the case where $T'_{vu}$  does not pass through $z$. By Claim~\ref{z} and Lemma~\ref{vertexnumber2+}, $T'_{vu}$ contains at most three vertices of $N_G(u)$. Suppose, for contradiction, that $T'_{vu}$ does not pass through $z$ and $|T_{vu} \cap N_G(u)| = 3$. Then, by Steps 11–16 of Algorithm~2, there exists exactly one vertex, say $z'$, between $T_{vu}$ and $H_X$, and $z' \in Z$.

    \item[(4)] In this case, the paths $T'_{su}$, $T_{vu}$, and the subgraph $H_X$ together contain at least six vertices of $d^*_G(u)$, contradicting the assumption that $d^*_G(u) \le 5$.
\end{itemize}

Therefore, if $T_{vu}$ does not pass through $z$, then $|T_{vu} \cap N_G(u)| \le 2$.

We now complete the argument by analyzing the paths from $u$ to each vertex $w \in V(D_u)$. We claim that $T_{uw}$ may contain three vertices of $N_G(u)$ if and only if it passes through $z$; otherwise, it contains at most two vertices of $N_G(u)$.

First, suppose the directed path $T'_{uw}$ in $\overrightarrow{H_u}$ does not contain $z$. Then $T'_{uw}$ remains unchanged in $D_u$ and satisfies $|T_{uw} \cap N_G(u)| \le 2$, so we are done in this case.

It remains to consider those vertices $w \in V(D_u)$ for which the corresponding path in $\overrightarrow{H_u}$ passes through $z$.  From the definition of $z$, there exists a vertex $y \in Y$ such that $Q_{yu}$ passes through $z$ and no other vertex of $Y$ except $y$. Then the subpath $uQ_{yu}(y, z)$ is a directed path from $u$ to $z$ in $D_u$, containing at most two vertices of $N_G(u)$. By construction, now we consider the two cases:
\begin{itemize}
    \item[(1)] If $w \in X$ and $P_{uw}$ contains $z$, then by Lemma~\ref{vertexnumber1+}, $P_{uw}$ contains at most two vertices of $N_G(u)$. Since $P_{uw}$ starts at $u$, we consider the subpath $P_{uw} - u$, which is a directed path from $z$ to $w$ in $D_u$, containing at most two vertices of $N_G(u)$. Therefore, the union of $uQ_{yu}(y, z)$ and $P_{uw} - u$ forms a path from $u$ to $w$ in $D_u$ containing at most three vertices of $N_G(u)$.
    \item[(2)] If $w \in Z$ and $Q_{uw}$ contains $z$, we distinguish two subcases:
    \begin{itemize}
        \item[(a)] If $w = z$, then $uQ_{yu}(y, z)$ is a valid path from $u$ to $w$ in $D_u$, containing at most two vertices of $N_G(u)$.
        \item[(b)] If $w \ne z$, then by Lemma~\ref{dvertexnumber3}, the original path $Q_{uw}$ in $\overrightarrow{H_u}$ contains at most two vertices of $N_G(u)$. Thus, the subpath $Q_{uw} - u$ forms a directed path from $z$ to $w$ in $D_u$, again containing at most two vertices of $N_G(u)$. The union of $uQ_{yu}(y, z)$ and $Q_{uw} - u$ gives a path from $u$ to $w$ in $D_u$ with at most three vertices of $N_G(u)$.
    \end{itemize}
\end{itemize}

Combining this with our previous result on $T_{vu}$, we conclude that for each pair $v, w \in V(D_u)$, there exist a $vu$-path $T_{vu}$ and a $uw$-path $T_{uw}$ such that
\[
|T_{vu} \cap N_G(u)| \le 3 \quad \text{and} \quad |T_{uw} \cap N_G(u)| \le 3,
\]
with $|T_{vu} \cap N_G(u)| = 3$ if and only if $z \in V(T_{vu})$, and $|T_{uw} \cap N_G(u)| = 3$ if and only if $z \in V(T_{uw})$.

Since both $T_{vu}$ and $T_{uw}$ contain $z \in N_G(u)$, their union forms a directed path from $v$ to $w$ in $D_u$ containing exactly five vertices of $N_G(u)$, which completes the case.

\noindent\textbf{Case 3.} $d^+_G(u) + d^-_G(u) \ge 2$ and $\min\{d^+_G(u), d^-_G(u)\} = 0.$

By symmetry, we may assume that $d^+_G(u) \ge 2$ and $d^-_G(u) = 0$. In this case, the neighborhood of $u$ is partitioned into $X_2$, $Y = Y_1 \cup Y_2$, and $Z$. For $v, w \in \overrightarrow{H_u}$, let $T'_{vw}$ denote a path from $v$ to $w$ in $\overrightarrow{H_u}$.

We first show that for each $v \in V(\overrightarrow{H_u})$, there exists a $vu$-path containing at most $d^+_G(u) + 3$ vertices of $N_G(u)$:
\begin{itemize}
    \item[(1)] If $v \in X$, then $T'_{vu} := vu$ contains exactly one vertex of $N_G(u)$;
    \item[(2)] If $v \in Y$, then $T'_{vu} := Q_{vu}$ contains at most $d^+_G(u) + 3$ vertices of $N_G(u)$ by Lemma~\ref{vertexnumber2+};
    \item[(3)] If $v \in Z$, then $T'_{vu} := Q_{vu}$ is a directed path from $v$ to $u$ in $\overrightarrow{H_u}$ that contains at most two vertices of $N_G(u)$ by Lemma~\ref{dvertexnumber3};
    \item[(4)] If $v \in V(\overrightarrow{H_u}) \setminus N_G[u]$, then there exists a vertex $v^\star \in N_G(u)$ such that the subpath of $T'_{v^\star u}$ from $v$ to $u$ serves as the desired path.
\end{itemize}

We next show that for each $w \in V(\overrightarrow{H_u})$, there exists a $uw$-path containing at most two vertices of $N_G(u)$:
\begin{itemize}
    \item[(1)] If $w \in X$, then $T'_{uw} := P_{uw}$ contains at most two vertices of $N_G(u)$ by Lemma~\ref{vertexnumber1+};
    \item[(2)] If $w \in Y$, then $T'_{uw} := uw$ contains exactly one vertex of $N_G(u)$;
    \item[(3)] If $w \in Z$, then $T'_{uw} := Q_{uw}$ is a directed path from $u$ to $w$ in $\overrightarrow{H_u}$ that contains at most two vertices of $N_G(u)$ by Lemma~\ref{dvertexnumber3};
    \item[(4)] If $w \in V(\overrightarrow{H_u}) \setminus N_G[u]$, then there exists a vertex $w^\star \in N_G(u)$ such that the subpath of $T'_{uw^\star}$ from $u$ to $w$ is the desired path.
\end{itemize}

If for every $v \in V(\overrightarrow{H_u})$, the path $T'_{vu}$ contains at most $d^+_G(u) + 2$ vertices of $N_G(u)$, then we may take $D_u := \overrightarrow{H_u}$, and we are done.

Suppose instead that there exists a vertex $s \in V(\overrightarrow{H_u})$ such that $T'_{su}$ contains exactly $d^+_G(u) + 3$ vertices of $N_G(u)$. The following claim shows that such a path must pass through some fixed vertex $z_0 \in Z$.

\begin{claim}\label{z_0}
There exists a vertex $z_0 \in Z$ such that for every $y \in Y_2$, if $Q_{yu}$ contains exactly $d^+_G(u) + 3$ vertices of $N_G(u)$, then $z_0$ is contained in $Q_{yu}$.
\end{claim}

\begin{proof}[Proof of Claim~\ref{z_0}]
Suppose there exists a vertex $s \in V(\overrightarrow{H_u})$ such that $T'_{su}$ contains exactly $d^+_G(u) + 3$ vertices of $N_G(u)$. By Steps 11–16 of the algorithm and Lemma~\ref{vertexnumber2}, we conclude that $s \in Y_2$, and the path $P_{su}$ contains all vertices of $Y_1$ and exactly one additional vertex $z_0 \in H_X$. Moreover, $z_0 \in Z$.

Now consider another vertex $t \in V(\overrightarrow{H_u})$ such that $T'_{tu}$ also contains exactly $d^+_G(u) + 3$ vertices of $N_G(u)$. Then, again by Steps 11–16 of the algorithm, Lemma~\ref{vertexnumber2}, and the definition of an economical path, we have $t \in Y_2$, and $P_{tu}$ contains all vertices of $Y_1$ and exactly one vertex $t' \in H_X$, with $t' \in Z$. Since $Y_1 \subseteq V(P_{su}) \cap V(P_{tu})$, Observation~\ref{shortestpath} implies that for each $y \in Y_1$, the subpath $P_{su}(y, u)$ is identical to $P_{tu}(y, u)$. Hence, $t' = z_0$.

Therefore, all such paths must pass through the same vertex $z_0 \in Z$.
\end{proof}

Let $D_u$ be the directed graph obtained from $\overrightarrow{H_u}$ by reversing the arc $\overrightarrow{uz_0}$ to $\overrightarrow{z_0u}$. We show that $D_u$ is our desired orientation.

We first show that for each $v \in V(D_u)$, there exists a directed path $T_{vu}$ from $v$ to $u$ containing at most $d^+_G(u) + 2$ vertices of $N_G(u)$. If the directed path $T'_{vu}$ in $\overrightarrow{H_u}$ does not pass through $z_0$, then $T'_{vu}$ is also a valid path in $D_u$, satisfying the desired bound. By Claim~\ref{z_0}, it suffices to consider $v = z_0$ and those vertices $v \in Y$ for which $Q_{vu}$ passes through $z_0$.

\begin{itemize}
    \item[(1)] For $v = z_0$, the path $T_{vu} := z_0u$ contains exactly one vertex of $N_G(u)$.
    \item[(2)] For $v \in Y$ such that $Q_{vu}$ contains $z_0$, we form a path by replacing the subpath $Q_{vu}(z_0, u)$ with the arc $z_0u$; this new path contains at most $d^+_G(u) + 2$ vertices of $N_G(u)$.
    \item[(3)] For $v \in V(D_u) \setminus N_G[u]$, there exists a vertex $v^\star \in N_G(u)$ such that the subpath of $T_{v^\star u}$ from $v$ to $u$ serves as the desired path.
\end{itemize}

Next, we show that for each $w \in V(D_u)$, there exists a directed path $T_{uw}$ from $u$ to $w$ containing at most three vertices of $N_G(u)$. If the path $T'_{uw}$ in $\overrightarrow{H_u}$ does not contain $z_0$, then it remains valid in $D_u$, and we are done.
We now consider the remaining cases: $w \in X$ such that $P_{uw}$ contains $z_0$, and $w \in Z$ such that $Q_{uw}$ contains $z_0$. By the definition of $z_0$, there exists a vertex $y \in Y$ such that $Q_{yu}$ includes $z_0$ and no other vertices of $Y$ besides $y$. Then the subpath $uQ_{yu}(y, z_0)$ is a path from $u$ to $z_0$ in $D_u$, containing at most two vertices of $N_G(u)$.

\begin{itemize}
    \item[(1)] For $w \in X$ such that $P_{uw}$ contains $z_0$, since $P_{uw}$ contains at most two vertices of $N_G(u)$ in $\overrightarrow{H_u}$ by Lemma~\ref{vertexnumber1+}, the subpath $P_{uw} - u$ is a directed path from $z_0$ to $w$ in $D_u$, containing at most two vertices of $N_G(u)$. The union of $uQ_{yu}(y, z_0)$ and $P_{uw} - u$ forms a path from $u$ to $w$ in $D_u$ containing at most three vertices of $N_G(u)$.

    \item[(2)] For $w \in Z$ such that $Q_{uw}$ contains $z_0$:
    \begin{itemize}
        \item[(a)] If $w = z_0$, then $uQ_{yu}(y, z_0)$ is the desired path, containing at most two vertices of $N_G(u)$.
        \item[(b)] If $w \ne z_0$, then since $Q_{uw}$ in $\overrightarrow{H_u}$ contains at most two vertices of $N_G(u)$ by Lemma~\ref{dvertexnumber3}, the subpath $Q_{uw} - u$ forms a path from $z_0$ to $w$ in $D_u$, also containing at most two vertices. Thus, the union of $uQ_{yu}(y, z_0)$ and $Q_{uw} - u$ is a path from $u$ to $w$ in $D_u$ containing at most three vertices of $N_G(u)$.
    \end{itemize}
    \item[(3)] For $w \in V(D_u) \setminus N_G[u]$, there exists a vertex $w^\star \in N_G(u)$ such that the subpath of $T_{uw^\star}$ from $u$ to $w$ is a valid path.
\end{itemize}

Therefore, for each $v, w \in V(D_u)$, we have paths $T_{vu}$ and $T_{uw}$ such that
\[
|T_{vu} \cap N_G(u)| \le d^+_G(u) + 2 \quad \text{and} \quad |T_{uw} \cap N_G(u)| \le 3.
\]

If $|T_{vu} \cap N_G(u)| + |T_{uw} \cap N_G(u)| \le d^+_G(u) + 4$, then $T_{vu} \cup T_{uw}$ is our desired path from $v$ to $w$, and we are done.

Otherwise, $v \in Y$ and $w \in X \cup Z$, where $T_{vu}$ contains exactly $d^+_G(u) + 2$ vertices and $T_{uw}$ contains exactly three vertices of $N_G(u)$. Suppose $Q_{vu}$ does not contain $z_0$. Then by Lemma~\ref{vertexnumber2+}, $T_{vu} := Q_{vu}$ contains at most three vertices of $N_G(u)$. But this contradicts $d^+_G(u) + 2 \ge 4 > 3$. Therefore, $T_{vu}$ must contain $z_0$.

If $T_{uw}$ does not contain $z_0$, then it coincides with $T'_{uw}$, which contains at most two vertices of $N_G(u)$—again a contradiction. Hence, both $T_{vu}$ and $T_{uw}$ must contain $z_0 \in N_G(u)$, and their union forms a path from $v$ to $w$ in $D_u$ that contains exactly $d^+_G(u) + 4$ vertices of $N_G(u)$. This completes the proof.

\noindent\textbf{Case 4.} $d^+_G(u) \ge 1$ and $d^-_G(u) \ge 1$.

In this case, the neighborhood of $u$ is partitioned into $X = X_1 \cup X_2$, $Y = Y_1 \cup Y_2$, and $Z$.

We first show that for each $v \in V(\overrightarrow{H_u})$, there exists a $vu$-path $T'_{vu}$ containing at most $d^+_G(u) + 3$ vertices of $N_G(u)$:
\begin{itemize}
    \item[(1)] If $v \in X$, then $T'_{vu} := vu$ contains exactly one vertex of $N_G(u)$;
    \item[(2)] If $v \in Y$, then $T'_{vu} := Q_{vu}$ contains at most $1 + d^+_G(u) + 2 = d^+_G(u) + 3$ vertices of $N_G(u)$ by Lemma~\ref{vertexnumber2+};
    \item[(3)] If $v \in Z$, then $T'_{vu} := Q_{vu}$ is a directed path from $v$ to $u$ in $\overrightarrow{H_u}$ that contains at most two vertices of $N_G(u)$ by Lemma~\ref{dvertexnumber3};
    \item[(4)] If $v \in V(\overrightarrow{H_u}) \setminus N_G[u]$, then there exists a vertex $v^\star \in N_G(u)$ such that the subpath of $T'_{v^\star u}$ from $v$ to $u$ is our desired path.
\end{itemize}

Similarly, for each $w \in V(\overrightarrow{H_u})$, there exists a $uw$-path $T'_{uw}$ containing at most $d^-_G(u) + 2$ vertices of $N_G(u)$:
\begin{itemize}
    \item[(1)] If $w \in X$, then $T'_{uw} := P_{uw}$ contains at most $d^-_G(u) + 2$ vertices of $N_G(u)$ by Lemma~\ref{vertexnumber1+};
    \item[(2)] If $w \in Y$, then $T'_{uw} := uw$ contains exactly one vertex of $N_G(u)$;
    \item[(3)] If $w \in Z$, then $T'_{uw} := Q_{uw}$ is a directed path from $u$ to $w$ in $\overrightarrow{H_u}$ that contains at most $d^-_G(u) + 2$ vertices of $N_G(u)$ by Lemma~\ref{dvertexnumber3};
    \item[(4)] If $w \in V(\overrightarrow{H_u}) \setminus N_G[u]$, then there exists a vertex $w^\star \in N_G(u)$ such that the subpath of $T'_{uw^\star}$ from $u$ to $w$ is our desired path.
\end{itemize}

If either of the following conditions holds:
\begin{itemize}
    \item[(1)] For every $v \in V(\overrightarrow{H_u})$, the path $T'_{vu}$ contains at most $d^+_G(u) + 2$ vertices of $N_G(u)$; or
    \item[(2)] For every $w \in V(\overrightarrow{H_u})$, the path $T'_{uw}$ contains at most $d^-_G(u) + 1$ vertices of $N_G(u)$,
\end{itemize}
then we may take $D_u := \overrightarrow{H_u}$ as the desired orientation, and we are done.

Suppose instead that there exist vertices $s, t \in V(\overrightarrow{H_u})$ such that $T'_{su}$ contains exactly $d^+_G(u) + 3$ vertices of $N_G(u)$ and $T'_{ut}$ contains exactly $d^-_G(u) + 2$ vertices of $N_G(u)$. By an analogous argument as in the previous case, we see that Claim~\ref{z_0} also holds. Moreover, the following claim is valid.
\begin{claim}\label{t_0}
There exists a vertex $t_0 \in Y \cup Z$ such that:
\begin{itemize}
    \item[\rm (1)] For every $x \in X_2$, if $P_{ux}$ contains exactly $d^-_G(u) + 2$ vertices of $N_G(u)$, then $t_0$ is contained in $P_{ux}$;
    \item[\rm (2)] For every $z \in Z$, if $Q_{uz}$ contains exactly $d^-_G(u) + 2$ vertices of $N_G(u)$, then $t_0$ is contained in $Q_{uz}$.
\end{itemize}
\end{claim}

\begin{proof}[Proof of Claim~\ref{t_0}]
Let $T'_{ut} := P_{ut}$ be a path that contains exactly $d^-_G(u) + 2$ vertices of $N_G(u)$. By Lemma~\ref{vertexnumber1}, $P_{ut}$ contains all vertices of $X_1$ and exactly one additional vertex, say $t_0$, from $Y \cup Z$.
Now consider another vertex $t' \in X_2$ such that $T'_{ut'}$ also contains $d^-_G(u) + 2$ vertices of $N_G(u)$. Again by Lemma~\ref{vertexnumber1}, the path $P_{ut'}$ contains all vertices of $X_1$ and exactly one vertex, say $t_1$, from $Y \cup Z$. Since $X_1 \subseteq V(P_{ut}) \cap V(P_{ut'})$, Observation~\ref{shortestpath} implies that the subpaths $P_{ut}(u,x)$ and $P_{ut'}(u,x)$ are identical for every $x \in X_1$. Hence, $t_1 = t_0$.

Next, consider a vertex $z \in Z$ such that $T'_{uz} := Q_{uz}$ contains exactly $d^-_G(u) + 2$ vertices of $N_G(u)$. By Observation~\ref{Q_{uz}}, we have $Q_{uz} := \overrightarrow{C^2_y}(u, z)$ where $C^2_y := P_{ux'}(u, v_y) P_{yu}(v_y, u)$, and $x' \in X$ is such that $P_{ux'}$ contains $v_y$ (the last vertex of $P_{yu}$ contained in $H_X$).
Since $P_{ux'}$ also contains all vertices of $X_1$ and $X_1 \subseteq V(P_{ut}) \cap V(P_{ux'})$, it follows from Observation~\ref{shortestpath} that the subpaths $P_{ut}(u,x)$ and $P_{ux'}(u,x)$ are identical for every $x \in X_1$. Therefore, the unique vertex $t_2 \in Y \cup Z$ in $Q_{uz}$ must also be $t_0$.

This completes the proof.
\end{proof}

If $z_0 = t_0$, we let $D_u := \overrightarrow{H_u}$ and define $T_{vw} := T'_{vw}$ for each $v, w \in V(D_u)$. If
\[
|V(T_{vu}) \cap N_G(u)| + |V(T_{uw}) \cap N_G(u)| \le d^+_G(u) + d^-_G(u) + 4,
\]
then $T_{vu} \cup T_{uw}$ is our desired path, and we are done. Otherwise, we have
\[
|V(T_{vu}) \cap N_G(u)| + |V(T_{uw}) \cap N_G(u)| = d^+_G(u) + d^-_G(u) + 5.
\]
By Claims~\ref{z_0} and~\ref{t_0}, this implies $v \in Y$, $w \in X \cup Z$, where $T_{vu}$ contains exactly $d^+_G(u) + 3$ vertices and $T_{uw}$ contains exactly $d^-_G(u) + 2$ vertices of $N_G(u)$. Since $z_0 = t_0$ and $z_0 \in V(T_{vu}) \cap V(T_{uw}) \cap N_G(u)$, the union $T_{vu} \cup T_{uw}$ is a directed path from $v$ to $w$ in $D_u$ containing $d^+_G(u) + d^-_G(u) + 4$ vertices of $N_G(u)$, and we are done.

If $z_0 \ne t_0$, let $D_u$ be the directed graph obtained from $\overrightarrow{H_u}$ by reversing the arc $\overrightarrow{uz_0}$ to $\overrightarrow{z_0u}$. We now show that $D_u$ is the desired orientation.

We claim that for every $v \in V(D_u)$, there exists a path $T_{vu}$ containing at most $d^+_G(u) + 2$ vertices of $N_G(u)$. If $T'_{vu}$ in $\overrightarrow{H_u}$ does not include $z_0$, it remains valid in $D_u$. By Claim~\ref{z_0}, we only need to consider $v = z_0$ and those $v \in Y$ for which $Q_{vu}$ contains $z_0$. For $v = z_0$, let $T_{vu} := z_0u$, which contains exactly one vertex of $N_G(u)$. For $v \in Y$ such that $Q_{vu}$ contains $z_0$, we replace the subpath $Q_{vu}(z_0,u)$ with the arc $z_0u$ to form a path that contains at most $d^+_G(u) + 2$ vertices of $N_G(u)$. For any $v \in V(D_u) \setminus N_G[u]$, there exists $v^\star \in N_G(u)$ such that the subpath of $T_{v^\star u}$ from $v$ to $u$ is the desired path.

Next, we show that for every $w \in V(D_u)$, there exists a path $T_{uw}$ containing at most $d^-_G(u) + 2$ vertices of $N_G(u)$. If $T'_{uw}$ does not contain $z_0$, then it is valid in $D_u$ as well. We need only consider $w \in X$ such that $P_{uw}$ contains $z_0$ and $w \in Z$ such that $Q_{uw}$ contains $z_0$. By definition of $z_0$, there is a vertex $y \in Y$ such that $Q_{yu}$ passes through $z_0$ and no other vertex in $Y \setminus \{y\}$. Since $z_0$ is the penultimate vertex of $Q_{yu}$, the subpath $uQ_{yu}(y, z_0)$ is a directed path from $u$ to $z_0$ in $D_u$ containing at most two vertices of $N_G(u)$.

For $w \in X$, if $P_{uw}$ contains a vertex in $X_1$, then by the construction process, it must contain $t_0$. Since $z_0 \ne t_0$ and $P_{uw}$ contains exactly one vertex in $Y \cup Z$, this leads to a contradiction. Thus, $P_{uw}$ contains no vertices in $X_1$, and by Lemma~\ref{vertexnumber1+}, it contains at most two vertices of $N_G(u)$. Therefore, $P_{uw} - u$ is a path from $z_0$ to $w$ in both $\overrightarrow{H_u}$ and $D_u$ with at most two vertices of $N_G(u)$. The union of $uQ_{yu}(y, z_0)$ and $P_{uw} - u$ forms a path in $D_u$ from $u$ to $w$ containing at most three vertices of $N_G(u)$.

For $w \in Z$ such that $Q_{uw}$ contains $z_0$, if $w = z_0$, then $uQ_{yu}(y, z_0)$ is our desired path. Otherwise, if $Q_{uw}$ contains a vertex of $X_1$, then by Observation~\ref{Q_{uz}}, we have $Q_{uw} = \overrightarrow{C^2_y}(u, w)$ where $C^2_y = P_{ux'}(u, v_y) P_{yu}(v_y, u)$ for some $x' \in X$ and $v_y \in H_X$. Then $P_{ux'}$ contains $t_0$, which again contradicts $z_0 \ne t_0$. Hence, $Q_{uw}$ contains no vertex in $X_1$, and by Lemma~\ref{dvertexnumber3}, $Q_{uw}$ contains at most two vertices of $N_G(u)$. The subpath $Q_{uw} - u$ is a path from $z_0$ to $w$ in both $\overrightarrow{H_u}$ and $D_u$, also with at most two vertices of $N_G(u)$. The union of $uQ_{yu}(y, z_0)$ and $Q_{uw} - u$ forms the desired $uw$-path in $D_u$.

For any $w \in V(D_u) \setminus N_G[u]$, there exists $w^\star \in N_G(u)$ such that the subpath of $T_{uw^\star}$ from $u$ to $w$ is valid. Hence, for every pair $v, w \in V(D_u)$, there exists a path $T_{vw}$ in $D_u$ from $v$ to $w$ containing at most $d^+_G(u) + d^-_G(u) + 4$ vertices of $N_G(u)$.

\medskip

\noindent This completes the proof of the theorem. \hfill$\qed$

We now construct several families of $(n, \Delta^*)$-graphs to demonstrate the tightness of the bounds established in Theorem~\ref{maximumdegree}. Dankelmann, Gao, and Surmacs~\cite{Dankelmann} introduced three families of undirected graphs to illustrate the sharpness of their results. We adopt their constructions in our setting.

\begin{definition}[\upshape Dankelmann, Gao, and Surmacs~\cite{Dankelmann}]
For $n \in \mathbb{N}$ with $n \ge 5$, let $\hat{O}_n$ be the graph obtained from the undirected cycle $u_1u_2 \dots u_nu_1$ by adding the edge $u_nu_2$. Let $S^0_n(4)$ be the graph obtained from the undirected cycle $uw_1 \dots w_{n-3}u$ and a path $P_2$ of order $2$ by joining $u$ to both vertices of $P_2$.
\end{definition}

\begin{definition}[\upshape Dankelmann, Gao, and Surmacs~\cite{Dankelmann}]
For $\Delta \ge 4$ and $n \ge 2\Delta + 2$, define $\psi := \frac{n}{2} - \Delta + 3$. The undirected equal-partition bigraph $R^{\Delta}_n = (V^{\Delta}_n, E^{\Delta}_n)$ is defined as follows:
\begin{align*}
V^{\Delta}_n &:= \{a_i, b_i : 1 \le i \le \frac{n}{2} \}, \\
E^{\Delta}_n &:= \{b_ia_i, a_ib_{i+1} : 1 \le i \le \tfrac{n}{2} - 1\} \cup \{b_1a_{\psi - 2}, a_{\psi - 2}b_\psi\} \cup \{b_\psi a_i : \psi \le i \le \tfrac{n}{2} \}.
\end{align*}
\end{definition}

The following remark and three theorems demonstrate that Theorem~\ref{maximumdegree} is optimal, except possibly in the case when a vertex $u$ with $d^*_G(u) = \Delta^*(G) = 5$ satisfies $d^+_G(u) + d^-_G(u) = 1$.

\begin{remark}
Sharp examples demonstrating the optimality of Theorem~\ref{maximumdegree} for undirected graphs have been constructed by Dankelmann, Gao, and Surmacs~\cite{Dankelmann}.
\end{remark}

For a positive integer $n$, let $P_n$ and $C_n$ denote the undirected $n$-path and $n$-cycle, respectively, and let $\overrightarrow{P_n}$ and $\overrightarrow{C_n}$ denote the directed $n$-path and $n$-cycle, respectively.

\begin{theorem}
For integers $n, \Delta^* \in \mathbb{N}$ with $0 \leq \Delta^* < n$, there exists an $(n, \Delta^*)$-graph with a vertex $u$ satisfying $d^*_G(u) = \Delta^*$ and $d^+_G(u) + d^-_G(u) = 0$, whose oriented diameter is at least $\min\{n-1, n - \Delta^* + 4\}$.
\end{theorem}

\begin{proof}
Since the desired mixed graph $G$ is bridgeless, we have $d_G(v) \ge 2$ for each vertex $v \in V(G)$. Given that $d^+_G(u) + d^-_G(u) = 0$, it follows that $d^*_G(u) = \Delta^* \ge 2$.
\begin{itemize}
    \item [(1)] For $\Delta^* = 2$, let $G_{n,2}$ be an undirected $n$-cycle.
    \item [(2)] For $\Delta^* = 3$, let $G_{n,3} := \hat{O}_n$.
    \item [(3)] For $\Delta^* = 4$, let $G_{n,4} := S^0_n(4)$.
    \item [(4)] For $\Delta^* = 5$, let $G_{n,5}$ be the graph obtained from the cycle $uu_1u_2\ldots u_{n-4}u$ and the path $P_3 = x_1x_2x_3$ by joining $u$ to each vertex of $P_3$, and orienting $x_ix_{i+1}$ from $x_i$ to $x_{i+1}$ for $i=1,2$.
\end{itemize}

It is easy to verify that each $G_{n,k}$ satisfies the required properties for $2 \le k \le 5$.

\smallskip

For $\Delta^* = 6$, let $G_{n,6}$ be the graph obtained from two cycles $u_1u_2u_3u_4u_1$ and $v_1v_2\ldots v_{n-3}v_1$ by identifying $u_1$ and $v_1$ into a new vertex $u$, adding edges $uu_3$ and $uv_3$, and orienting: $u_i u_{i+1}$ from $u_i$ to $u_{i+1}$ for $i=2,3$;  $v_j v_{j+1}$ from $v_j$ to $v_{j+1}$ for $2 \le j \le n-4$.
In any strong orientation $D_{n,6}$ of $G_{n,6}$, we have
\[
\overrightarrow{uu_2},\; \overrightarrow{u_4u},\; \overrightarrow{uv_2},\; \overrightarrow{v_{n-3}u} \in E(D_{n,6}).
\]
Then, depending on the directions of the remaining edges incident to $u$:
\begin{itemize}
    \item [(1)] If $\overrightarrow{uu_3}, \overrightarrow{uv_3} \in E(D_{n,6})$, then $d_{D_{n,6}}(u_2, v_{n-3}) = n - 2$;
    \item [(2)] If $\overrightarrow{uu_3}, \overrightarrow{v_3u} \in E(D_{n,6})$, then $d_{D_{n,6}}(u_2, v_{n-3}) = n - 1$;
    \item [(3)] If $\overrightarrow{u_3u}, \overrightarrow{uv_3} \in E(D_{n,6})$, then $d_{D_{n,6}}(v_2, u_4) = n - 1$;
    \item [(4)] If $\overrightarrow{u_3u}, \overrightarrow{v_3u} \in E(D_{n,6})$, then $d_{D_{n,6}}(u_2, v_{n-3}) = n - 2$.
\end{itemize}
Thus, $\overrightarrow{\mathrm{diam}}(G_{n,6}) \ge n - 2 = n - \Delta^* + 4$.

\begin{center}
\scalebox{0.8}[0.8]{\includegraphics{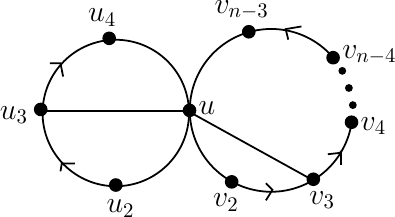}}\\
\small Fig. 1. The graph $G_{n,6}$.
\end{center}

\smallskip

For $\Delta^* \ge 7$, construct $G_{n,\Delta^*}$ as follows: Start with an undirected cycle $uu_1\ldots u_{n-\Delta^*+1}u$; Add a directed path $\overrightarrow{P_3} = x_1x_2x_3$; Add some copies of $P_2$ and $P_3$, joined to $u$, so that the total number of vertices in these paths equals $\Delta^* - 5$.

\begin{center}
\scalebox{0.8}[0.8]{\includegraphics{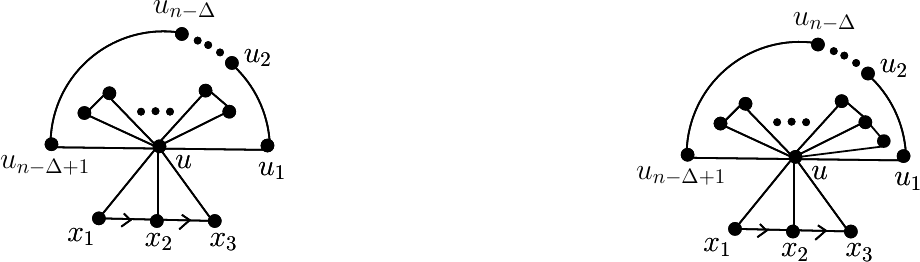}}\\
\small Fig. 2. The graph $G_{n,\Delta^*}$ with $\Delta^* \ge 7$.
\end{center}

In any strong orientation $D_{n,\Delta^*}$ of $G_{n,\Delta^*}$, we have:
\[
\overrightarrow{ux_1}, \overrightarrow{x_3u}, \overrightarrow{uu_1}, \overrightarrow{u_{n - \Delta^* + 1}u}, \overrightarrow{u_iu_{i+1}} \in E(D_{n,\Delta^*}) \quad \text{for } 1 \le i \le n - \Delta^*.
\]
Then:
\begin{itemize}
    \item [(1)]  If $\overrightarrow{x_2u} \in E(D_{n,\Delta^*})$, then $d_{D_{n,\Delta^*}}(u_1, x_3) = n - \Delta^* + 4$;
    \item [(2)] If $\overrightarrow{ux_2} \in E(D_{n,\Delta^*})$, then $d_{D_{n,\Delta^*}}(x_1, u_{n - \Delta^* + 1}) = n - \Delta^* + 4$.
\end{itemize}

Therefore, $\overrightarrow{\mathrm{diam}}(G_{n,\Delta^*}) \ge n - \Delta^* + 4$.
\end{proof}

\begin{theorem}
For integers $n, \Delta^* \in \mathbb{N}$ with $0 \leq \Delta^* < n$, there exists an $(n, \Delta^*)$-graph with a vertex $u$ satisfying $d^*_G(u) = \Delta^* \ne 5$ and $d^+_G(u) + d^-_G(u) = 1$, whose oriented diameter is at least $\min\{n-1, n-\Delta^*+4\}$.
\end{theorem}

\begin{proof}
Since our desired mixed graph $G$ is bridgeless, we have $d_G(v) \ge 2$ for each vertex $v \in V(G)$. It follows from $d^+_G(u) + d^-_G(u) = 1$ that $d^*_G(u) = \Delta^* \ge 1$.

We construct examples for small values of $\Delta^*$ as follows:
\begin{itemize}
    \item[(1)] For $\Delta^* = 1$, let $F_{n,1}$ be the oriented $n$-cycle obtained by directing all edges of an $n$-cycle $u_1u_2\ldots u_nu_1$ as $\overrightarrow{u_i u_{i+1}}$ for $1 \le i \le n-1$.

    \item[(2)] For $\Delta^* = 2$, let $F_{n,2}$ be the graph obtained from $\hat{O}_n$ by orienting the edge $u_2u_n$ from $u_2$ to $u_n$.

    \item[(3)] For $\Delta^* = 3$, let $F_{n,3}$ be the graph obtained from $S^0_n(4)$ by orienting $uu_1$ from $u$ to $u_1$.

    \item[(4)] For $\Delta^* = 4$, let $F_{n,4}$ be the graph formed from the cycle $uu_1\ldots u_{n-4}u$ and the path $P_3 = x_1x_2x_3$ by joining $u$ to every vertex of $P_3$, orienting $uu_1$ from $u$ to $u_1$ and $x_i x_{i+1}$ from $x_i$ to $x_{i+1}$ for $i = 1, 2$.
\end{itemize}

For $\Delta^* = 5$, let $F_{n,5}$ be the graph obtained from $G_{n,6}$ by orienting $uu_3$ from $u$ to $u_3$. Then $\overrightarrow{\mathrm{diam}}(F_{n,5}) \ge n - 2 = n - \Delta^* + 3$.

\begin{center}
\scalebox{0.8}[0.8]{\includegraphics{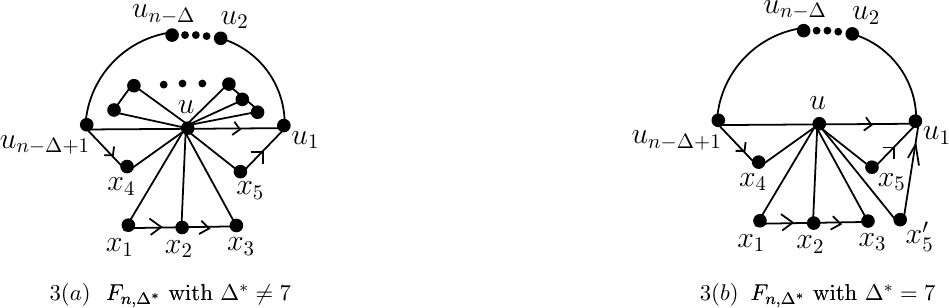}}\\
\small Fig 3. The graph $F_{n,\Delta^*}$ with $\Delta^* \ge 6$.
\end{center}

Now consider $\Delta^* \ge 6$, with $\Delta^* \ne 7$. Let $J_{\Delta^*}$ be the graph formed from: an undirected cycle $uu_1\ldots u_{n-\Delta^*+1}u$,  a directed path $\overrightarrow{P_3} = x_1x_2x_3$,
a stable set $\{x_4, x_5\}$ joined to $u$,
 additional edges $u_{n-\Delta^*+1}x_4$ (oriented from $u_{n-\Delta^*+1}$ to $x_4$) and $x_5u_1$ (oriented from $x_5$ to $u_1$),
 the edge $uu_1$ oriented from $u$ to $u_1$.

Let $F_{n,\Delta^*}$ be the graph obtained by joining $u$ to the vertices of disjoint copies of $P_2$ and $P_3$ paths such that the total number of added vertices is $\Delta^* - 6$. The case $\Delta^* = 7$ is depicted in Fig. 3(b).

In any strong orientation $D_{n,\Delta^*}$ of $F_{n,\Delta^*}$, we have $\overrightarrow{ux_1}, \overrightarrow{x_3u}, \overrightarrow{x_4u}, \overrightarrow{ux_5} \in E(D_{n,\Delta^*})$ and  $\overrightarrow{u_iu_{i+1}} \in E(D_{n,\Delta^*})$ for $1 \le i \le n - \Delta^*$.

Consider the four possible cases of orientations involving $u_{n - \Delta^* + 1}$ and $x_2$:
\begin{itemize}
    \item[(1)] If $\overrightarrow{u_{n - \Delta^* + 1} u}, \overrightarrow{u x_2} \in E(D_{n,\Delta^*})$, then $d_{D_{n,\Delta^*}}(x_1, x_4) = n - \Delta^* + 4$.
    \item[(2)] If $\overrightarrow{u_{n - \Delta^* + 1} u}, \overrightarrow{x_2 u} \in E(D_{n,\Delta^*})$, then $d_{D_{n,\Delta^*}}(x_5, x_3) = n - \Delta^* + 4$.
    \item[(3)] If $\overrightarrow{u u_{n - \Delta^* + 1}}, \overrightarrow{u x_2} \in E(D_{n,\Delta^*})$, then $d_{D_{n,\Delta^*}}(x_5, x_3) = n - \Delta^* + 4$.
    \item[(4)] If $\overrightarrow{u u_{n - \Delta^* + 1}}, \overrightarrow{x_2 u} \in E(D_{n,\Delta^*})$, then $d_{D_{n,\Delta^*}}(x_1, x_4) = n - \Delta^* + 5$.
\end{itemize}

In all cases, the oriented diameter satisfies $\overrightarrow{\mathrm{diam}}(F_{n,\Delta^*}) \ge n - \Delta^* + 4$.
\end{proof}

\begin{theorem}
For integers $n, \Delta^* \in \mathbb{N}$ with $0 \leq \Delta^* < n$, there exists an $(n, \Delta^*)$-graph with a vertex $u$ satisfying $d^*_G(u) = \Delta^*$ and $d^+_G(u) + d^-_G(u) \geq 2$, whose oriented diameter is at least $\min\{n-1, n - \Delta^* + 3\}$.
\end{theorem}

\begin{proof}
We construct examples based on the value of $\Delta^*$:

\begin{itemize}
    \item[(1)] For $\Delta^* = 0$, let $M_{n,0}$ be a directed cycle of order $n$.

    \item[(2)] For $\Delta^* = 1$, let $M_{n,1}$ be the graph obtained from the path $x_1x_2\ldots x_n$ by orienting each edge from $x_i$ to $x_{i+1}$ for $1 \le i \le n-1$.

    \item[(3)] For $\Delta^* = 2$, let $M_{n,2}$ be the graph obtained from $\hat{O}_n$ by orienting the edges $u_i u_{i+1}$ from $u_i$ to $u_{i+1}$, for $1 \le i \le n-1$.

    \item[(4)] For $\Delta^* = 3$, let $M_{n,3}$ be the graph obtained from the bipartite graph $R^{\Delta}_n$ by orienting each $a_i b_i$ from $a_i$ to $b_i$, and $b_{\psi} a_j$ from $b_{\psi}$ to $a_j$, where $1 \le i \le \psi - 1$ and $\psi + 1 \le j \le \frac{n}{2}$.
\end{itemize}

Each $M_{n,k}$ constructed above satisfies the required properties for $0 \le k \le 3$.

\begin{center}
\scalebox{0.8}[0.8]{\includegraphics{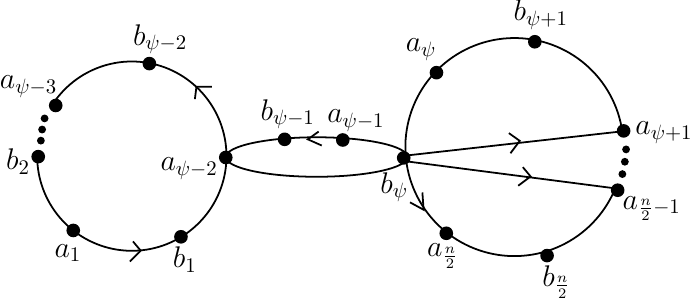}}\\
\small Fig 4. $M_{n,3}$.
\end{center}

\begin{center}
\scalebox{0.8}[0.8]{\includegraphics{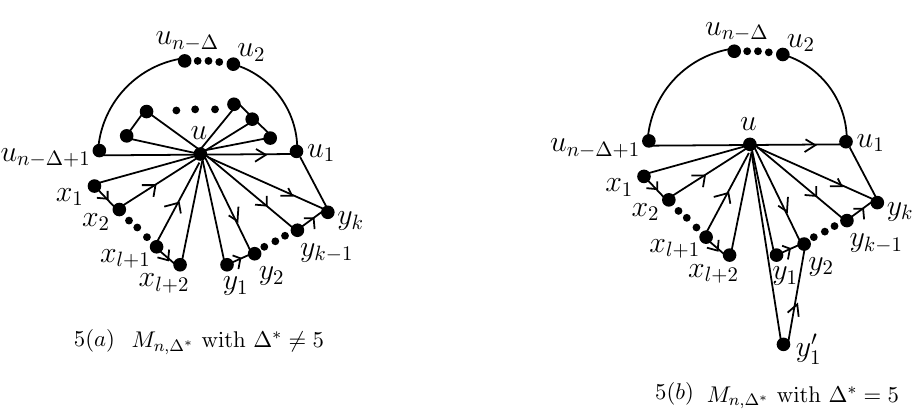}}\\
\small Fig 5. The graph $M_{n,\Delta^*}$ with $\Delta^* \ge 4$.
\end{center}

For $\Delta^* \ge 4$ and $\Delta^* \ne 5$, we define $M_{n,\Delta^*}$ as follows:

\begin{itemize}
    \item[(1)] Let $S_{\Delta^*}$ be the graph consisting of an undirected cycle $uu_1\ldots u_{n-\Delta^*+1}u$,
 a directed path $\overrightarrow{P_{l+2}} = x_1x_2\ldots x_{l+2}$,
 a directed path $\overrightarrow{P_k} = y_1y_2\ldots y_k$,
 edges $ux_i$ (oriented from $x_i$ to $u$) for $2 \le i \le l+1$,
 edges $uy_j$ (oriented from $u$ to $y_j$) for $2 \le j \le k$,
 the edge $u_1y_k$ (oriented from $y_k$ to $u_1$),
 the edge $uu_1$ (oriented from $u$ to $u_1$).

    \item[(2)] Let $M_{n,\Delta^*}$ be the graph obtained by taking $S_{\Delta^*}$ and attaching some $P_2$ and $P_3$ paths such that $u$ is joined to every vertex of these paths, and the total number of added vertices is $\Delta^* - 5$. Here, $l = d^-_{M_{n,\Delta^*}}(u)$ and $k = d^+_{M_{n,\Delta^*}}(u)$.

    \item[(3)] For $\Delta^* = 5$, the construction is given in Fig. 4(b).
\end{itemize}

In any strong orientation $D_{n,\Delta^*}$ of $M_{n,\Delta^*}$, the following arcs must appear:
 $\overrightarrow{ux_1}$, $\overrightarrow{uy_1}$, $\overrightarrow{x_{l+2}u}$, and $\overrightarrow{y_k u_1}$,
 $\overrightarrow{u_{i} u_{i+1}}$ for $1 \le i \le n - \Delta^*$,
$\overrightarrow{u_{n - \Delta^* + 1} u}$.

Then the path from $y_1$ to $x_{l+2}$ passes through $u_1, u_2, \ldots, u_{n - \Delta^* + 1}, u, x_{l+2}$, resulting in a path of length $n - \Delta^* + 3$. Thus,
\[
\overrightarrow{\mathrm{diam}}(M_{n,\Delta^*}) \ge n - \Delta^* + 3.
\]

This completes the proof.
\end{proof}

\subsection{The proof of Theorem \ref{bipartite}}

In this subsection, we prove Theorem \ref{bipartite} and present a family of sharp examples for the case $d^+_G(u) + d^-_G(u) \ge 2$.

\begin{proof}[Proof of Theorem \ref{bipartite}.]
We consider two main cases based on whether $G$ is undirected or a general mixed graph.

\begin{itemize}
    \item[(1)] \text{$G$ is an undirected bipartite graph.} \\
    By Theorem~\ref{Dbipartite}, $G$ has an our desired strongly orientation.

    \item[(2)] \text{$G$ is a bridgeless mixed bipartite graph.} \\
    By the same reasoning as in the proof of Theorem~2.6, there exists a strong orientation of $G$ such that for every fixed $u \in B$:
    \begin{itemize}
        \item[(a)] If $d^+_G(u) + d^-_G(u) \ge 2$, then any two vertices can be connected by a path containing at most $d^+_G(u) + d^-_G(u) + 4$ vertices of $N_G(u)$.
        \item[(b)] If $d^+_G(u) + d^-_G(u) < 2$, then such a path contains at most $d^+_G(u) + d^-_G(u) + 5$ vertices of $N_G(u)$.
    \end{itemize}
    Since each such path must alternate between $A$ and $B$, it contains at most $(|A| - d^*_G(u) + 4)$ vertices in  $A$, and at most $(|A| - d^*_G(u) + 5)$ vertices in $B$, for a total of at most $2(|A| - \Delta^*) + 9$ vertices when $d^+_G(u) + d^-_G(u) \ge 2$, and $2(|A| - \Delta^*) + 11$ otherwise.
\end{itemize}

This concludes the proof.
\end{proof}

We now show that the bound in Theorem~\ref{bipartite} is tight, except in the special case when $G$ is a mixed bipartite graph in which each vertex $u$ with maximum undirected degree satisfies $d^+_G(u) + d^-_G(u) \le 1$.

\begin{remark}
Sharp examples demonstrating the tightness of Theorem~\ref{bipartite} in the case of undirected bipartite graphs have been constructed by Dankelmann, Gao, and Surmacs~\cite{Dankelmann}.
\end{remark}

\begin{theorem} \label{R}
For integers $n, \Delta^* \in \mathbb{N}$ with $5 \leq \Delta^* < n$, there exists a mixed equal-partition bipartite graph of order $n$ with a vertex $u$ satisfying $d^*_G(u) = \Delta^*$, whose oriented diameter is at least $n - 2\Delta^* + 8$.
\end{theorem}

\begin{proof}
Let $l, k \in \mathbb{N}$ with $l + k \le \frac{n}{2} - 2$ and $r \le \frac{n}{2} - 3$.

 Define $H_k$ as follows: start with the cycles $x_1x_2\ldots x_{2k+4}x_1$ and $a_0a_1a_2a_3a_0$, identify $x_2$ and $a_0$ into a new vertex $u_2$, add edges $x_1x_{2+2d}$ for $1 \le d \le k$, and orient $u_2x_3$ and $x_ix_{i+1}$ from $x_i$ to $x_{i+1}$ for $3 \le i \le 2k+3$.

 Define $H_l$ similarly from the cycles $y_1y_2\ldots y_{2l+4}y_1$ and $b_0b_1b_2b_3b_0$, identify $y_{2l+4}$ and $b_0$ into a vertex $b_4$, add edges $y_1y_{2+2d}$ for $1 \le d \le l$, and orient $y_{2l+3}b_4$ and $y_iy_{i+1}$ from $y_i$ to $y_{i+1}$ for $2 \le i \le 2l+2$.

 Let $H_r$ be the undirected graph obtained from the cycle $z_1z_2\ldots z_{2r-8}z_1$ by adding edges $z_1z_{2+2d}$ for $1 \le d \le r - 6$.

Now construct $R_{n, \Delta^*}$:
\begin{itemize}
    \item [(1)] If $\Delta^* = 5$, let $G = R_{n,5}$ be formed from $H_k$ and $H_l$ by identifying $x_1$ and $y_1$ into a new vertex $u_1$, and add a vertex $a_4$ with edges $a_4u_1$ and $a_4x_3$, where $a_4x_3$ is oriented from $a_4$ to $x_3$.
    \item [(2)] If $\Delta^* \ge 6$, let $G = R_{n,\Delta^*}$ be obtained by identifying $x_1$, $y_1$, and $z_1$ into a new vertex $u_1$ from the disjoint union of $H_k$, $H_l$, and $H_r$, with $k = d^+_G(u_1)$, $l = d^-_G(u_1)$, $r = d^*_G(u_1)$, and $n = 2(k + l + r + 2)$.
\end{itemize}

\begin{center}
\scalebox{0.7}[0.7]{\includegraphics{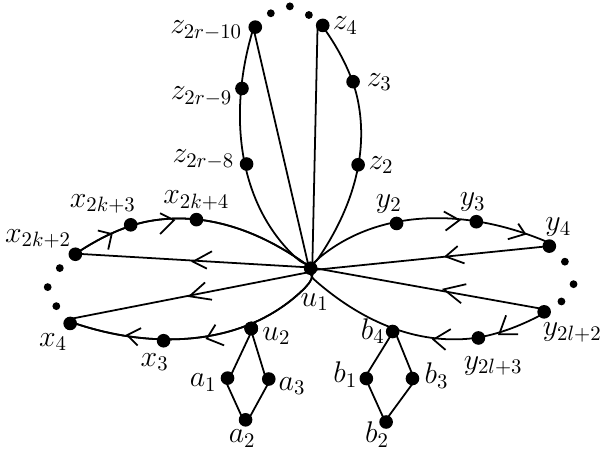}}

\small Fig 6. The construction of $R_{n,\Delta^*}$.
\end{center}

Note that $R_{n,\Delta^*}$ contains no odd cycles and hence is bipartite. Assign all vertices with even indices to partite set $A$, and all with odd indices to partite set $B$. Since each addition preserves parity, we have $|A| = |B|$, making $G$ a mixed equal-partition bipartite graph.

In any strong orientation $D_{n,\Delta^*}$ of $G$, one has
$\overrightarrow{u_1x_2}, \overrightarrow{x_{2k+4}u_1}, \overrightarrow{u_1y_2}, \overrightarrow{y_{2l+4}u_1} \in E(D_{n,\Delta^*})$.
By symmetry, with loss of any generality, assume that $\overrightarrow{u_2a_1}, \overrightarrow{a_3u_2}, \overrightarrow{b_4b_1} \in E(D_{n,\Delta^*})$ and $\overrightarrow{a_ia_{i+1}}, \overrightarrow{b_jb_{j+1}} \in E(D_{n,\Delta^*})$ for $1 \le i \le 2$ and $1 \le j \le 3$.

Then, the directed distance from $a_1$ to $b_3$ in $D_{n,\Delta^*}$ is $2(k + l) + 12 = n - 2\Delta^* + 8$, and we conclude:
\[
\overrightarrow{\mathrm{diam}}(R_{n,\Delta^*}) \ge n - 2\Delta^* + 8.
\]
This completes the proof.
\end{proof}

\begin{remark}
By Theorem~\ref{R}, we conclude that Theorem~\ref{bipartite} is optimal, except possibly when $G$ is a bridgeless mixed equal-partition bipartite graph containing at least one directed edge, and every vertex $u$ with the maximum undirected degree satisfies $d_G(u) \leq d^*_G(u) + 1$.
\end{remark}

\section{Conclusion}

In this paper, we investigated the oriented diameter of bridgeless mixed graphs with bounded undirected degree. For a mixed graph \( G \) and a vertex \( u \in V(G) \) with maximum undirected degree \( d^*_G(u) = \Delta^*(G) \), we provided explicit constructions and bounds on the length of a directed path connecting any two vertices in a strong orientation of \( G \).

Our main result shows that there exists a strong orientation \( D \) of \( G \) such that for any pair of vertices \( v, w \in V(G) \), there exists a \( vw \)-path \( P_{vw} \) in \( D \) with
\[
|V(P_{vw}) \cap N_G(u)| \le f(\Delta^*(G)),
\]
where the bound \( f(\Delta^*(G)) \) depends on the sum \( d^+_G(u) + d^-_G(u) \). We classified the cases based on this sum and gave tight bounds:
\begin{itemize}
    \item[(1)] If \( d^+_G(u) + d^-_G(u) = 0 \), then \( f(\Delta^*) \le \min\{n - 1, n - \Delta^* + 4\} \);
    \item[(2)] If \( d^+_G(u) + d^-_G(u) = 1 \), then \( f(\Delta^*) \le \min\{n - 1, n - \Delta^* + 4\} \);
    \item[(3)] If \( d^+_G(u) + d^-_G(u) \ge 2 \), or \( d^+_G(u) + d^-_G(u) = 1 \) and $d^*(u)\le 5$, then \( f(\Delta^*) \le \min\{n - 1, n - \Delta^* + 3\} \).
\end{itemize}

We also studied mixed bipartite graphs with equal partition \( V(G) = A \cup B \), and showed that for any fixed \( u \in B \), there exists a strong orientation \( D \) such that any pair of vertices is connected by a directed path \( P \) with
\[
|V(P) \cap N_G(u)| \le \max\{2(|A| - d^*_G(u)) + 9, \, 2(|A| - d^*_G(u)) + 11\},
\]
depending on whether \( d^+_G(u) + d^-_G(u) \ge 2 \) or not. Furthermore, we constructed families of mixed graphs to show that these upper bounds are optimal in all but one special case: when \( d^+_G(u) + d^-_G(u) \le 1 \).

Our results open several natural questions for further study:

\begin{itemize}
 \item[(1)]  We study the problem of minimizing the oriented diameter for bridgeless mixed graphs of a given maximum undirected degree. Thanks to our constructive proof, it is not difficult to derive a polynomial-time algorithm. Whether there exists a polynomial-time $c$-approximation algorithm, where $c$ is a constant, is an interesting question for further research. 
    
 \item[(2)] Extend the analysis to weighted mixed graphs: For an edge-weighted bridgeless mixed graph, study the upper bound of its oriented diameter in terms of the maximum undirected degree and edge weights.
    
 \item[(3)] For general graphs (not necessarily bridgeless), how does the presence of bridges or cut-edges affect the oriented diameter relative to \( \Delta^*(G) \)?
   That is, only non-bridge edges are oriented, while bridge edges remain undirected, aiming to minimize the diameter.

\end{itemize}

\section{Acknowledgements}
This work was supported by Natural Science Foundation of Henan Province (No. 252300421488), the Postgraduate Education Reform and Quality Improvement Project of Henan Province (No.YJS2024KC32), and the Young Backbone Teachers of Higher Education Institutions in Henan Province  (No. 2024GGJS045).


\begin{thebibliography}{25}

\bibitem{Bosak} J. Bos\'{a}k, A. Rosa, S. Zn\'{a}m. On decompositions of complete graphs into factors with given diameters. In: \emph{Theory of Graphs (Proc. Colloq., Tihany, 1966)}, Academic Press, New York, pp. 37--56, 1968.

\bibitem{Boesch} F. Boesch, R. Tindel. Robbins's theorem for mixed multigraphs. \emph{Amer. Math. Monthly} \textbf{87}(10): 716--719, 1980.

\bibitem{Bondy} J. A. Bondy, U. S. R. Murty. \emph{Graph Theory}. GTM 244, Springer, 2008.

\bibitem{Bau} S. Bau, P. Dankelmann. Diameter of orientations of graphs with given minimum degree. \emph{European J. Combin.} \textbf{49}: 126--133, 2015.

\bibitem{Babu} J. Babu, D. Benson, D. Rajendraprasad. Improved bounds for the oriented radius of mixed multigraphs. \emph{J. Graph Theory} \textbf{103}(4): 674--689, 2023.

\bibitem{Chvatal} V. Chv\'{a}tal, C. Thomassen. Distances in orientations of graphs. \emph{J. Combin. Theory Ser. B} \textbf{24}(1): 61--75, 1978.

\bibitem{Chung} F. R. K. Chung, M. R. Garey, R. E. Tarjan. Strongly connected orientations of mixed multigraphs. \emph{Networks} \textbf{15}(4): 477--484, 1985.

\bibitem{Czabarka} \'{E}. Czabarka, P. Dankelmann, L. A. Sz\'{e}kely. A degree condition for diameter two orientability of graphs. \emph{Discrete Math.} \textbf{342}(4): 1063--1065, 2019.

\bibitem{Chen} B. Chen, A. Chang. Oriented diameter of graphs with given girth and maximum degree. \emph{Discrete Math.} \textbf{346}(4): 113287, 2023.

\bibitem{Cochran} G. Cochran. Large girth and small oriented diameter graphs. \emph{Discrete Math.} \textbf{347}(4): 113846, 2024.

\bibitem{Dankelmann} P. Dankelmann, Y. Guo, M. Surmacs. Oriented diameter of graphs with given maximum degree. \emph{J. Graph Theory} \textbf{88}(1): 5--17, 2018.

\bibitem{P.Dankelmann} P. Dankelmann, J. Morgan, E. Rivett-Carnac. The oriented diameter of graphs with given connected domination number and distance domination number. \emph{Graphs Combin.} \textbf{40}(1): 18, 2024.

\bibitem{Fomin} F. V. Fomin, M. Matamala, E. Prisner, I. Rapaport. AT-free graphs: Linear bounds for the oriented diameter. \emph{Discrete Appl. Math.} \textbf{141}(1): 135--148, 2004.

\bibitem{Huang} X. Huang, H. Li, X. Li, Y. Sun. Oriented diameter and rainbow connection number of a graph. \emph{Discrete Math. Theor. Comput. Sci.} \textbf{16}(3): 51--60, 2014.

\bibitem{Kwok} P. K. Kwok, Q. Liu, D. B. West. Oriented diameter of graphs with diameter 3. \emph{J. Combin. Theory Ser. B} \textbf{100}(3): 265--274, 2010.

\bibitem{Kurz} S. Kurz, M. L\"{a}tsch. Bounds for the minimum oriented diameter. \emph{Discrete Math. Theor. Comput. Sci.} \textbf{14}(1): 109--140, 2012.

\bibitem{Kumar} K. S. A. Kumar, D. Rajendraprasad, K. S. Sudeep. Oriented diameter of star graphs. \emph{Discrete Appl. Math.} \textbf{319}: 362--371, 2022.

\bibitem{Li} H. Li, Z. Ding, J. Liu, Y. Gao, S. Zhao. A note on improved bounds for the oriented radius of mixed multigraphs. arXiv:2407.01612 [math.CO], 2024.

\bibitem{Mondal} D. Mondal, N. Parthiban, I. Rajasingh. On the oriented diameter of planar triangulations. \emph{J. Combin. Optim.} \textbf{47}(5): 79, 2024.

\bibitem{Robbins} H. E. Robbins. A theorem on graphs with an application to a problem of traffic control. \emph{Amer. Math. Monthly} \textbf{46}(5): 281--283, 1939.

\bibitem{Surmacs} M. Surmacs. Improved bound on the oriented diameter of graphs with given minimum degree. \emph{European J. Combin.} \textbf{59}: 187--191, 2017.

\bibitem{Wang1} X. Wang, Y. Chen, P. Dankelmann, Y. Guo, M. Surmacs, L. Volkmann. Oriented diameter of maximal outerplanar graphs. \emph{J. Graph Theory} \textbf{98}(3): 426--444, 2021.

\bibitem{Wang2} X. Wang, Y. Chen. Optimal oriented diameter of graphs with diameter 3. \emph{J. Combin. Theory Ser. B} \textbf{155}: 374--388, 2022.

\end{thebibliography}
\end{document}